\documentclass[preprint,review,11pt]{elsarticle}

\usepackage{amssymb}
\usepackage{amsmath}
\usepackage{amsthm}
\usepackage{esint}
\newtheorem{theorem}{Theorem}[section]

\numberwithin{equation}{section}

\newtheorem{lemma}[theorem]{Lemma}

\theoremstyle{definition}
\newtheorem{remark}[theorem]{Remark}

\journal{}

\usepackage{geometry}

\begin{document}

\baselineskip15pt

\begin{frontmatter}



\title{\huge Interior Second Order H\"{o}lder Regularity for Stokes systems\tnoteref{t1}}
\tnotetext[t1]{Research supported by National Natural Science Foundation of China [grant 12001419 and grant 12071365], Fundamental Research Funds for the Central Universities [grant xzy012022004] and Shaanxi Fundamental Science Research Project for Mathematics and Physics
[grant 22JSZ003 and grant 23JSQ026].\\
\indent \mbox{}\hskip0.3cm Data sharing is not applicable to this article as no new data were created or analyzed in this study.\\
\indent \mbox{}\hskip0.3cm All authors declare that they have no potential conflict of interest.
}


\address[rvt1]{School of Mathematics and Statistics, Xi'an Jiaotong University, Xi'an 710049, China}
\address[rvt2]{Department of Mathematics, University of Iowa, Iowa City, IA52242, USA}

\author[rvt1]{Rong Dong}
\ead{dongrong1203@mail.xjtu.edu.cn }
\author[rvt1]{Dongsheng Li\corref{cor1}}
\ead{lidsh@mail.xjtu.edu.cn}
\cortext[cor1]{Corresponding author.}
\author[rvt2]{Lihe Wang}
\ead{lihe-wang@uiowa.edu}

\begin{abstract}
 Global second order H\"{o}lder regularity for Stokes systems can be obtained by
global Schauder estimates, which are actually a priori estimates and were established by Solonnikov \cite{s11} and \cite{s16} with appropriate compatible conditions. This paper will investigate the corresponding interior regularity which unfortunately may fail in general from Serrin's counterexample (cf. \cite{s14}). However, we discover interior $C^{2,\alpha}$ regularity for velocity and interior $C^{1,\alpha}$ regularity for pressure in spatial variables, and furthermore,
for curl of velocity, we find its gradient belongs to $C^{\alpha, \frac{\alpha}2}$,
that is, possesses H\"{o}lder continuity in both space and time directions.
The interesting phenomenon here is that
 no continuity in time variable is assumed for both the coefficients and the righthand side terms.
The estimates for velocity and its curl are achieved pointwisely and
the results are sharp indicated by a counterexample.
\end{abstract}

\bigskip


\end{frontmatter}


\section{Introduction}

The primary purpose of this paper is to investigate interior second order H\"{o}lder regularity for Stokes systems with variable coefficients in non-divergence form. The system concerned here is precisely the following
\begin{equation}\label{s0.0}
\left\{
\begin{array}{ll}
\vspace{2mm} u_t+\mathcal{L} u+\nabla p=f~~\mbox{in}~~\Omega_{T};\\
\vspace{2mm}
\mbox{div}~u=g~~\mbox{in}~~\Omega_{T},
\end{array}
\right.
\end{equation}
 where
\begin{equation*}
\mathcal{L}:=-A(x,t)D^2=:-a^{ij}(x,t)D_{ij},
\end{equation*}
$\Omega_{T}=\Omega\times(-T,0]$ is a cylinder with $T>0$ and $\Omega\subset\mathbb{R}^{n}$ ($n\geq2$), the velocity field $u(x,t)=(u_1, u_2,..., u_n)$ and the pressure $p(x,t)$ are unknowns, and $f$ and $g$ are respectively a vector field and a function given suitably.
Throughout this paper,  the summation convention is always assumed and the coefficient matrix $A(x,t)=\big(a^{ij}(x,t)\big)_{n\times n}$ satisfies the hypothesis

\noindent
($\bf H$) Ellipticity: $A\in L_{\infty}(\Omega_T)$ such that
\begin{equation*}
|a^{ij}(x,t)|\leq\lambda^{-1}~\mbox{and}~ a^{ij}(x,t)\xi_i\xi_j\geq\lambda|\xi|^2~ a.e. ~\mbox{in} ~\Omega_T
\end{equation*}
for some $\lambda\in(0,1]$ and any $\xi=(\xi_i)\in\mathbb{R}^{n}$.


Stokes systems in the form (1.1)  arise from non-Newton fluid (c.f.  e.g. \cite{BP}) and also from Stokes systems on manifolds (c.f.  e.g. \cite{DM} and \cite{MT}). Both (1.1) and its corresponding divergent form are recently studied in \cite{DKP}-\cite{DL5}.

Schauder estimates play an important role in the study of Stokes systems and e.g. from it, Abe and Giga \cite{s1} deduced the analyticity of the Stokes semigroup in $L^\infty_{\sigma}(\Omega)$.
For Stokes systems with constant coefficients, Solonnikov \cite{s11} established the global
Schauder estimates by using potential theory, and later, he further developed the
estimates to systems with variable coefficients (cf. \cite{s16} and the references therein), where
appropriate compatible conditions and suitable smoothness are assumed.
Chang and Kang \cite{CK} gave a global H\"{o}lder estimate for the velocity and then they can solve the system in H\"{o}lder spaces for constant coefficients and velocity with vanishing divergence.

From global Schauder estimates, Solonnikov obtained the existence of solutions of class $C^{2+\alpha,1+\frac\alpha2}$
and as an important consequence, the global second order H\"{o}lder regularity can be established.
Since Stokes systems are parabolic with degeneration (cf. \cite{sj1}), the global regularity heavily depends on the boundary and compatible conditions (e.g. seeing the global results mentioned above)
and the interior regularity may fail generally. The latter can be seen from
Serrin's counterexample (cf. \cite{s14} or the following Remark 1.2 (ii)). The reason is that,
roughly speaking, Stokes systems contain no regularity in the potential vector fields
in time direction.
In this paper, we will try to investigate the possible interior second order H\"{o}lder regularity and
actually, we will establish the interior $C^{2,\alpha}$ regularity for velocity and
$C^{1,\alpha}$ regularity for pressure in the spatial variables.
An interesting result is that for curl of velocity, the regularity is obtained in both spatial variables and time variable, that is, the $C^{\alpha,\frac{\alpha}2}$ regularity for  $D\nabla\times u$, while we only assume that
$A$, $f$ and $Dg$ are H\"{o}lder continuous in spatial variables.
The counterexamples mentioned above indicate that our results are sharp.
Comparing to our results,  Naumann and Wolf \cite{NW} derived an interior H\"{o}lder estimate for velocity in spatial direction for (1.1) with constant coefficients and vanishing divergence of velocity.

As far as to the authors' knowledge, there has been almost no interior second order H\"{o}lder regularity
for (1.1) until now. Only interior $L_p$ regularity is investigated, where the solution $u\in L_p$ is always assumed to compensate the loss of regularity in time direction for Stokes systems and this assumption is essential by Serrin's counterexample (cf. \cite{s14} and the discussion in \cite{HLW5}).
 Hu, Li and Wang \cite{HLW5} establish interior $L_p$ regularity of $D^2 u$ for
 (\ref{s0.0}) with constant coefficients and vanishing divergence of velocity, where if $u\in L_p$ is not assumed,
  only $\Delta u\in L_p$ can be obtained. Later, Dong, Li and Wang \cite{DL5} generalize the results for H\"{o}lder continuous coefficients in spatial variables.
 Dong and Phan \cite{DP5} obtain an a priori interior $L_p$ mixed-norm estimate, where the coefficients are of small mean oscillations in the spatial variables and recently, Dong and Kwon \cite{DK5} further achieve interior $L_p$ mixed-norm regularity under the same assumption.  Dong, Kim and Phan \cite{DKP} give $L_p$ mixed-norm estimates up to the boundary for (1.1) with Lions type boundary conditions, where since the estimates are not global, the suitable integrability of $u$ is also assumed.

\bigskip
The following are the main results of this paper.

\bigskip

\begin{theorem}
Suppose $M_0>0$, $0<\alpha<1$,
and $A\in C_{x}^{\alpha}(\Omega_{T})$ with
$[A]_{x,\alpha;\Omega_T}\leq M_0$ and $(H)$ holding.
Let
$u\in W^{2,1}_2(\Omega_T)$,
$Dp\in L_2(\Omega_T)$, $f\in C_{x}^{\alpha}(\Omega_{T})$ and $g\in C_{x}^{1,\alpha}(\Omega_{T})$ satisfy (\ref{s0.0}).
For any
$0<T'<T$ and any $\Omega'\subset\subset\Omega$, we have

(i)
$D^2u\in C_x^{\alpha}(\overline{\Omega'_{T'}})$ and $D\nabla\times u\in C^{\alpha,\frac{\alpha}{2}}(\overline{\Omega'_{T'}})$
 such that
\begin{eqnarray}\label{g01}
\begin{aligned}
&\|D^2 u\|_{C_x^{\alpha}(\overline{\Omega'_{T'}})}+\|D\nabla\times u\|_{C^{\alpha,\frac{\alpha}{2}}(\overline{\Omega'_{T'}})}\\
&\mbox{}\hskip0.2cm\leq
C\Big(\|u\|_{L_\infty(-T,0;L_2(\Omega))}+\|D\nabla\times u\|_{L_2(\Omega_{T})}
+\|f\|_{C_{x}^{\alpha}(\Omega_{T})}
+\|g\|_{C_{x}^{1,\alpha}(\Omega_{T})}\Big);
\end{aligned}
\end{eqnarray}

(ii) If moreover, $p\in L_\infty(-T,0;L_2(\Omega))$ and $g_t\in L_\infty(-T,0;L_{q}(\Omega))$ with $q=\frac{n}{1-\alpha}$, then $Dp\in C_x^{\alpha}(\overline{\Omega'_{T'}})$ such that
\begin{eqnarray}\label{g02}
\begin{aligned}
&\|Dp\|_{C_x^{\alpha}(\overline{\Omega'_{T'}})}
\leq
C\Big(\|u\|_{L_\infty(-T,0;L_2(\Omega))}+\|D\nabla\times u\|_{L_2(\Omega_{T})}
+\|f\|_{C_{x}^{\alpha}(\Omega_{T})}\\
&\mbox{}\hskip2.5cm+\|g\|_{C_{x}^{1,\alpha}(\Omega_{T})}
+\|g_t\|_{L_\infty(-T,0;L_{q}(\Omega))}+\|p\|_{L_\infty(-T,0;L_2(\Omega))}\Big),
\end{aligned}
\end{eqnarray}
where $C$ depends only on $n$, $\lambda$, $\alpha$, $M_0$ and $\min\{\mbox{dist}(\Omega',\partial\Omega),\sqrt{T-T'}\}$.
\end{theorem}

\bigskip

\begin{remark}

(i) The existence of solutions in the theorem can be obtained globally by
 \cite{sj2}, where except for appropriate boundary and compatible conditions, the coefficient matrix $A$ needs further satisfying
the continuity with respect to $t\in[-T,0]$ and belonging to $W^{1}_q(\Omega)$ with $q>\frac{n}{2}$ for all fixed $t\in[-T,0]$.

(ii) Recall the counterexample given by Serrin \cite{s14}. Let $h$ be a harmonic function in $\Omega$ and $\sigma$ be a function of $L_1(-T,0)$. Then $u(x,t)=\nabla h(x)\int_{-T}^t\sigma(s)ds$ and $p(x,t)=-h(x)\sigma(t)$ satisfy (1.1) with $A=I$ the identity matrix, and vanishing $f$ and $g$.
It is easy to see that $D^2u=D^3 h(x)\int_{-T}^t\sigma(s)ds$ may not be H\"{o}lder continuous with respect to variable $t$. That is, our conclusions are sharp.

(iii) The pointwise version of conclusion (i) can be achieved (see Theorem 3.1 in Section 3) and only $u\in W^{2,1}_2(\Omega_T)$ is needed ($Dp\in L_2(\Omega_T)$ is not needed), which implies $u\in L_\infty(-T,0;L_2(\Omega))$ (cf. \cite{E1}).

(iv) Observe that $D\nabla\times u$ has H\"{o}lder continuity in both $x$ and $t$ variables although $A$, $f$ and $g$ are not assumed to be continuous in $t$.

(v) The condition $g_t\in L_\infty(-T,0;L_{q}(\Omega))$ in Theorem 1.1 (ii) can be also replaced by $g_t\in L_\infty(-T,0;C^{-1,\alpha}(\Omega))$, that is, $g_t=\mbox{div}~G$ for some $G\in L_{\infty}(-T,0; C^{\alpha}(\Omega))$.


\end{remark}

\bigskip

After Theorem 1.1 (i) is established, (ii) follows by solving the pressure from the Stokes system. As for (1.2), we will actually establish a pointwise version. That is, we will approximate in each scale the velocity and
its curl's gradient by a second order polynomial of spatial variables and a constant tensor
   in an error of $(2+\alpha)$-th order and of $\alpha$-th order of the scale respectively, where coefficients of the polynomials may depend on time.
  The approximation is indeed Campanato's characterization of $C^{k,\alpha}$ continuity (cf. \cite{s10} \cite{GG} \cite{L} and the following Lemma 2.1) and the polynomials and the constant tensors will be obtained by iteration.

  This method originates from Caffarelli \cite{s8} (cf. also \cite{s0}), however some essential modification is needed here. Since the polynomials are not of constant coefficients and may not be differentiable with respect to the time variable, we can not deal with the difference between the solution and the polynomial in each scale as it is done in \cite{s8}, where the difference satisfies a suitable equation. We can only rescale the equation and the solution and then decompose the solution several times according to the structure of the equation
  to utilize appropriately assumptions derived from  the previous scale.
  Though the sequence of polynomials are obtained by iteration,  the one of the new scale is not derived from the one of the previos scale by readjusting as in \cite{s8}.

 In the iteration to deduce the pointwise version of (1.2), the regularity is gained since the curl of velocity satisfies parabolic equations and the gradient of
  velocity can be solved from its curl and divergence by elliptic equations.
  This kind of manipulation for getting regularity of Stokes systems is also adopted for the interior $L_p$ estimates of (\ref{s0.0}) (cf. \cite{DK5} \cite{DP5} and \cite{DL5}).

The remainder of the paper will be organized as follows. Section 2 is concerned with some preliminary tools and the existence of solutions of (\ref{s0.0}) with the coefficient matrix $A$ depending only on the time $t$. In Section 3, we give the pointwise second order H\"{o}lder regularity of the velocity and H\"{o}lder regularity of the gradient of the curl. The proof of Theorem 1.1 (ii) will be presented in Section 4.

\bigskip

Before the end of this section, we introduce some notations.
\vskip 2mm\noindent{\bf Notations:}

$B_r(x):=\{y\in\mathbb{R}^n:|x-y|<r\}$;
~$B_r:=B_r(0)$.
\vspace{2mm}

$Q_r(x,t):=B_r(x)\times(t-r^2,t]$;~
$Q_r:=B_r\times(-r^2,0]$.
\vspace{2mm}

$\overline{(u)}_{Q_r(x,t)}:=\fint_{Q_r(x,t)}u(y,s)
dyds=\frac{1}{|Q_r(x,t)|}\int_{Q_r(x,t)}u(y,s)
dyds$.
\vspace{2mm}

$\overline{(u)}_{B_r(x)}(s):=\fint_{B_r(x)}u(y,s)
dy=\frac{1}{|B_r(x)|}\int_{B_r(x)}u(y,s)dy$.
\vspace{2mm}

$D_iu:=\frac{\partial u}{\partial x_i}$; $Du:=\big(D_{1}u,D_{2}u,...,D_{n}u\big)$; $D^2u:=\big(D_iD_ju\big)_{n\times n}$.
\vspace{2mm}

$\nabla\times u:=\big(D_ju_i-D_iu_j\big)_{n\times n}$.
\vspace{2mm}

$\partial_p \Omega_T: \mbox{parabolic~boundary~of}~\Omega_T$.
\vspace{2mm}

$[u]_{\alpha,\frac{\alpha}{2};\Omega_T}:=\sup_{(x,t),(y,s)\in \Omega_T, (x,t)\neq (y,s)}\frac{|u(x,t)-u(y,s)|}{|x-y|^{\alpha}+|t-s|^{\frac{\alpha}{2}}}$.
\vspace{2mm}

$\|u\|_{C^{\alpha,\frac{\alpha}{2}}(\Omega_T)}:=[u]_{\alpha,\frac{\alpha}{2};\Omega_T}+\|u\|_{L_\infty(\Omega_T)}$.

\vspace{2mm}

$C^{\alpha,\frac{\alpha}{2}}(\Omega_T):=\{u:\|u\|_{C^{\alpha,\frac{\alpha}{2}}(\Omega_T)}<\infty\}$.
\vspace{2mm}

$[u]_{x,\alpha;\Omega_T}:=\sup_{(x,t),(y,t)\in \Omega_T, x\neq y}\frac{|u(x,t)-u(y,t)|}{|x-y|^{\alpha}}$.
\vspace{2mm}

$\|u\|_{C_{x}^\alpha(\Omega_T)}:=[u]_{x,\alpha;\Omega_T}+\|u\|_{L_\infty(\Omega_T)}$.

\vspace{2mm}

$C^\alpha_{x}(\Omega_T):=\{u:\|u\|_{C_{x}^\alpha(\Omega_T)}<\infty\}$.
\vspace{2mm}

$\|u\|_{C_{x}^{1,\alpha}(\Omega_T)}:=[Du]_{x,\alpha;\Omega_T}+\|u\|_{L_\infty(\Omega_T)}+\|Du\|_{L_\infty(\Omega_T)}$.
\vspace{2mm}

$C^{1,\alpha}_{x}(\Omega_T):=\{u:\|u\|_{C_{x}^{1,\alpha}(\Omega_T)}<\infty\}$.

\section{Preliminary}

We first introduce the lemma below, which characterizes the H\"{o}lder continuity by Campanato spaces and we refer to \cite{L} and  \cite{s10} for its proof.

\bigskip

\begin{lemma}
Let $T>0$, $0<\alpha<1$ and $\Omega\subset\mathbb{R}^n$ be a bounded Lipschitz domain.

(i) If
\begin{eqnarray*}
\sup_{(x,t)\in\Omega_{T}}\sup_{r>0}
\Big(\fint_{Q_r(x,t)\cap\Omega_{T}}|u(y,s)-\overline{(u)}_{Q_r(x,t)\cap\Omega_{T}}|^2dyds\Big)^{\frac{1}{2}}\leq
Ar^{\alpha}
\end{eqnarray*}
for some constant $A$,
then $u\in C^{\alpha,\frac{\alpha}{2}}(\overline{\Omega_{T}})$ and
$$[u]_{\alpha,\frac{\alpha}{2};\overline{\Omega_{T}}}\leq C
A$$
for some constant $C$ depending only on $\alpha$, $n$ and $\Omega_T$.

(ii) If there exist constants $r_0, A>0$ and functions $a(x,t)$, $b(x,t)$ and $c(x,t)$ bounded in $\Omega_T$ such that
\begin{eqnarray}\label{j2.1}
\sup_{0<r\leq r_0}
\Big(\fint_{B_r(x)\cap\Omega}|u(y,t)-a(x,t)-b(x,t)y-\frac12c(x,t)y^2|^2dy\Big)^{\frac{1}{2}}\leq
Ar^{2+\alpha}
\end{eqnarray}
for any $(x,t)\in\Omega_T$,
then $D^2u\in C_x^{\alpha}(\overline{\Omega_{T}})$ and
$$[D^2u]_{x,\alpha;\overline{\Omega_{T}}}\leq C
A$$
for some constant $C$ depending only on $\alpha$, $n$, $T$ and $\Omega$.
\end{lemma}

\bigskip

\begin{remark} Lemma 2.1 (i) is for anisotropic spaces while (ii) is essentially not since the integral is only taking with respect to spatial variables and  $t$ is fixed in (\ref{j2.1}).
We write $\frac12cy^2$ instead of $\frac12c^{ij}y_iy_j$ for simplicity in (\ref{j2.1}) and from now on, we always do so.	
\end{remark}

\bigskip

The next Lemma is actually an equivalence between norms of polynomials, which may be well known.

\bigskip

\begin{lemma}
Let $T>0$, $a(t)\in L_{\infty}(-T,0;\mathbb{R})$, $b(t)\in L_{\infty}(-T,0;\mathbb{R}^n)$ and $c(t)\in L_{\infty}(-T,0;\mathbb{R}^{n^2})$.
If there exist $A$ and $r>0$ such that
\begin{equation}\label{j06}
\begin{aligned}
\sup_{t\in(-T,0]}\left(\fint_{B_{r}}|a(t)+b(t)x+\frac12c(t)x^2|^2dx\right)^{\frac12}
\leq A,
\end{aligned}
\end{equation}
then
\begin{equation}\label{j07}
\begin{aligned}
\sup_{t\in(-T,0]}|a(t)|
\leq CA,~~\sup_{t\in(-T,0]}|b(t)|
\leq C\frac{A}{r}~~\mbox{and}~~\sup_{t\in(-T,0]}|c(t)|
\leq C\frac{A}{r^2}
\end{aligned}
\end{equation}
for some constant $C$ depending only on $n$.
\end{lemma}

\vskip 2mm\noindent{\bf Proof.}\quad
For any second order polynomial $a+bx+\frac{1}{2}cx^2$,
there exists a constant $C$ depending only on $n$ such that
$$C^{-1}\left(\int_{B_1}|a+bx+\frac{1}{2}cx^2|^2dx\right)^{\frac12}\leq\max\{|a|,|b|,|c|\}\leq C\left(\int_{B_1}|a+bx+\frac{1}{2}cx^2|^2dx\right)^{\frac12}$$
or
$$C^{-1}\left(\fint_{B_r}|a+bx+\frac{1}{2}cx^2|^2dx\right)^{\frac12}\leq\max\{|a|,r|b|,r^2|c|\}\leq C\left(\fint_{B_r}|a+bx+\frac{1}{2}cx^2|^2dx\right)^{\frac12}$$
for any $r>0$. Then (\ref{j07}) is a consequence of (\ref{j06}).
\qed

\bigskip

The following two lemmas are $L_p$ estimates for elliptic and parabolic equations respectively and we include them here for convenience. The latter one follows easily from \cite{K} (see also \cite{DP5} and \cite{DL5}).

\bigskip

\begin{lemma}
Let $1<p<\infty$, $k\geq0$  and $\tilde\Omega\subset\subset\Omega$.
If $u\in W^{k+2}_{p}(\Omega;\mathbb{R}^{n})$, then
\begin{eqnarray*}
\begin{aligned}
\|D^{k+2}u\|_{L_p(\tilde{\Omega})}\leq&\ C\Big(\|D^k\mbox{div}~(\nabla\times u)\|_{L_p(\Omega)}
+\|D^{k+1}(\mbox{div}~u)\|_{L_p(\Omega)}+\|u\|_{L_1(\Omega)}\Big),
\end{aligned}
\end{eqnarray*}
where $C$ depends only on $k$, $n$, $p$, $\tilde\Omega$ and $\Omega$.
\end{lemma}
\vskip 2mm\noindent{\bf Proof.}\quad
Recall $\nabla\times u=D_ju_i-D_iu_j$ and then
\begin{eqnarray*}
\begin{aligned}
\Delta u
=\mbox{div}~\nabla\times u+\nabla~\mbox{div}~u
~~\mbox{in}~~\Omega.
\end{aligned}
\end{eqnarray*}
We obtain the lemma immediately by standard $L_p$ estimates for elliptic equations.
\qed

\bigskip

\begin{lemma}
Let $2\leq p<\infty$, $0<T'< T$ and $\tilde\Omega\subset\subset\Omega$. Suppose $(a^{ij}(t))_{n\times n}$ (not depending on $x$) satisfies $(H)$ and $u\in L_2(-T,0;H^1(\Omega))\cap L_{\infty}(-T,0;L_2(\Omega))$ is a weak solution of the following equation:
\begin{equation}\label{j10}
u_t-a^{ij}(t)D_{ij}u=g~~\mbox{in}~~\Omega_T.
\end{equation}

(i) If $g=0$, then $u, Du, u_t\in L_{\infty}(\tilde\Omega_{T'})$ and
\begin{eqnarray*}
\|u\|_{L_\infty(\tilde{\Omega}_{T'})}+
\|Du\|_{L_\infty(\tilde{\Omega}_{T'})}+\|u_t\|_{L_\infty(\tilde{\Omega}_{T'})}\leq C\|u\|_{L_p(\Omega_{T})},
\end{eqnarray*}
where $C$ depends only on $n$, $\lambda$, $p$, $T'$, $T$, $\tilde\Omega$ and $\Omega$.

(ii) If $g=\mbox{div}~f$ for some $f\in L_{p}(\Omega_{T};\mathbb{R}^{n})$, then $Du\in L_{p}(\tilde\Omega_{T'})$ and
\begin{eqnarray*}
\|Du\|_{L_p(\tilde{\Omega}_{T'})}\leq C\Big(\|f\|_{L_p(\Omega_{T})}+\|u\|_{L_p(\Omega_{T})}\Big),
\end{eqnarray*}
where $C$ depends only on $n$, $\lambda$, $p$, $T'$, $T$, $\tilde\Omega$ and $\Omega$.
\end{lemma}

\vskip 2mm\noindent{\bf Proof.}\quad
(ii) is an easy consequence of Theorem 2.4 in \cite{K} while (i) follows from (ii) by differentiating the equation. Actually, take $\varepsilon>0$ and let $\rho_\varepsilon$ be the scaled mollifier in $\mathbb{R}^{n}$. Let $\Omega^*\subset\Omega$ with smooth boundary, $dist(\Omega^*,\partial \Omega)>\varepsilon$ and $\tilde\Omega_{T'}\subset\subset\Omega^*_T$. Then by $g=0$,
\begin{equation*}
(u_{\varepsilon})_t-a^{ij}(t)D_{ij}u_\varepsilon=0~~\mbox{in}~~\Omega^*_{T},
\end{equation*}
where $u_\varepsilon=\rho_\varepsilon\ast u$.
Differentiating the equation suitable times with respect to $x$ and by (ii), we have
\begin{eqnarray*}
\|u_\varepsilon\|_{L_\infty(\tilde{\Omega}_{T'})}+\|Du_\varepsilon\|_{L_\infty(\tilde{\Omega}_{T'})}+\|u_{\varepsilon t}\|_{L_\infty(\tilde{\Omega}_{T'})}\leq C\|u_\varepsilon\|_{L_p(\Omega^*_{T})},
\end{eqnarray*}
where $C$ depends only on $n$, $\lambda$, $p$, $\tilde\Omega_{T'}$ and $\Omega_{T}$.
Letting $\varepsilon\rightarrow0$, (i) is obtained.
\qed

\bigskip

The last lemma concerns the solvability of the initial boundary-value problem for Stokes system with coefficient $A$ depending only on $t$. Since we can not find suitable references, a complete proof is given.

\begin{lemma}
Let $\Omega$ be a bounded Lipschitz domain. Suppose $A(t)=(a^{ij}(t))_{n\times n}$ (not depending on $x$) satisfies $(H)$ and $f\in L_2(\Omega_{T})$. The following Stokes system:
\begin{eqnarray}\label{jh0}
\left\{
\begin{array}{ll}
u_t-A(t)D^2u+\nabla p=f,~~\mbox{div}~u=0~~\mbox{in}~~\Omega_{T};\\
\vspace{2mm}
u=0~~\mbox{on}~~\left(\partial\Omega\times[-T,0]\right)\cup\left(\Omega\times\{-T\}\right)
\end{array}
\right.
\end{eqnarray}
has a weak solution
 $u\in L_2(-T,0;H_0^1(\Omega))\cap C([-T,0];L_2(\Omega))$ and $p=\frac{\partial}{\partial t}P$ for some $P\in C([-T,0];L_2(\Omega))$, where $u$ is unique
such that
\begin{eqnarray}\label{jh1}
\begin{aligned}
\sup_{t\in(-T,0]}\|u\|_{L_2(\Omega)}+\|u\|_{L_2(-T,0;H_0^1(\Omega))}
\leq C\|f\|_{L_2(\Omega_{T})}
\end{aligned}
\end{eqnarray}
for some constant $C$ depending only on $n$, $\lambda$, $\Omega$ and $T$.
\end{lemma}

\begin{remark}
Here $f$ can be released to be of class $L_2(-T,0;H^{-1}(\Omega))$.
\end{remark}

\vskip 2mm\noindent{\bf Proof of Lemma 2.6.}\quad
We adopt the arguments for proving Theorem III.1.1 and Proposition III.1.1 in \cite{tj1}.
First by the standard Galerkin's method, we obtain
$u\in L_2(-T,0;H_0^1(\Omega))$ with $u'\in L_2(-T,0;H^{-1}(\Omega))$ such that
\begin{eqnarray}\label{xw7}
<u',v>-<a^{ij}(t)D_{ij}u,v>=<f,v>~a.e.~\mbox{in}~(-T,0)
\end{eqnarray}
for any $v\in H_0^1(\Omega)$ with $\mbox{div}~v=0$, where $<\cdot,\cdot>$ denotes the dual product between $H^{-1}(\Omega)$ and $H_0^1(\Omega)$. Then we have $u\in C([-T,0];L_2(\Omega))$. To obtain the pressure $p$, we take integral on $[-T,t]$ for any $t\in[-T,0]$ on both side of (\ref{xw7}). It follows from the vanishing initial value of $u$ that
\begin{eqnarray*}
\Big<u(t)-\int_{-T}^t\big(a^{ij}(s)D_{ij}u+f\big)ds,v\Big>=0.
\end{eqnarray*}
Since $v\in H_0^1(\Omega)$ with $\mbox{div}~v=0$ is arbitrary, we have, by Proposition I.1.1 and I.1.2   in \cite{tj1}, that for each $t\in[-T,0]$, there exists $P(t)\in L_2(\Omega)$ such that
\begin{eqnarray*}
u(t)-\int_{-T}^t\big(a^{ij}(s)D_{ij}u+f\big)ds=-\nabla P(t).
\end{eqnarray*}
We deduce $\nabla P(t)\in C([-T,0];H^{-1}(\Omega))$ and then $P(t)\in C([-T,0];L_2(\Omega))$ by the continuity of the lefthand side with $t$.
\qed

\section{Pointwise estimates of velocity}\setcounter{equation}{0}

In this section, we will mainly show the following pointwise regularity estimates from which
Theorem 1.1 (i) follows clearly.
\begin{theorem}
Let $M_0,M_1,M_2>0$, $0<\alpha<1$ and $q>\frac{n}{1-\alpha}$. Suppose $Q_R(x_0,t_0)\subset\Omega_{T}$, $(H)$ holds and
$u\in W^{2,1}_2(\Omega_T)$,
$Dp\in L_2(\Omega_T)$, $f\in L_2(\Omega_T)$ and $g\in L_\infty(-T,0;W^{1}_q(\Omega))$ satisfying (\ref{s0.0}). If for any $0<r\leq R$,
\begin{equation}\label{t1}
\big|A(x,t)-A(x_0,t)\big|\leq M_0r^{\alpha}~~\mbox{in}~~Q_{r}(x_0,t_0),
\end{equation}
\begin{equation}\label{t2}
\fint_{Q_{r}(x_0,t_0)}|f-\overline{(f)}_{B_{r}(x_0)}|^2dxdt\leq M_1r^{2\alpha}
\end{equation}
and
\begin{equation}\label{t33}
\sup_{t\in(t_0-r^{2},t_0]}\fint_{B_{r}(x_0)}|Dg-\overline{(Dg)}_{B_{r}(x_0)}|^qdx\leq M_2r^{q\alpha},
\end{equation}
then there exist two constants $0<\sigma<1$ and $C>0$ depending only on $n$, $\lambda$, $\alpha$, $q$ and $M_0$, and three functions $a(t)\in L_\infty(t_0-\sigma^2 R^2,t_0;\mathbb{R}^{n})$, $b(t)\in L_\infty(t_0-\sigma^2 R^2,t_0;\mathbb{R}^{n^2})$ and $c(t)\in L_\infty(t_0-\sigma^2 R^2,t_0;\mathbb{R}^{n^3})$ such that
\begin{equation*}
\begin{aligned}
||a||_{L_\infty(t_0-\sigma^2 R^2,t_0)}+
||b||_{L_\infty(t_0-\sigma^2 R^2,t_0)}+
||c||_{L_\infty(t_0-\sigma^2 R^2,t_0)}\leq
 C{\cal K},
 \end{aligned}
\end{equation*}
\begin{equation}\label{p01}
\begin{aligned}
\fint_{Q_{r}(x_0,t_0)}|D\nabla\times u-\overline{(D\nabla\times u)}_{Q_{r}(x_0,t_0)}|^2dxdt
\leq Cr^{2\alpha}{\cal K}^2
\end{aligned}
\end{equation}
and
\begin{equation}\label{p02}
\begin{aligned}
\sup_{t\in(t_0-r^{2},t_0]}\fint_{B_{r}(x_0)}|u-a(t)-b(t)x-\frac12c(t)x^2|^2dx
\leq&\ Cr^{2(2+\alpha)}{\cal K}^2
\end{aligned}
\end{equation}
for any $0<r\leq\sigma R$, where
\begin{equation*}
\begin{aligned}
{\cal K}=&\|u\|_{L_\infty(t_0-R^2,t_0;L_2(B_R(x_0)))}+\|D\nabla\times u\|_{L_2(Q_R(x_0,t_0))}\\
&+\|g\|_{L_\infty(t_0-R^2,t_0;W^{1}_q(B_R(x_0)))}+M^\frac{1}{2}_1+M^\frac{1}{q}_2.
\end{aligned}
\end{equation*}
\end{theorem}

\bigskip
\begin{remark}
Observe that $u\in W^{2,1}_2(\Omega_T)$ implies
$\|u\|_{L_\infty(t_0-R^2,t_0;L_2(B_R(x_0)))}+\|D\nabla\times u\|_{L_2(Q_R(x_0,t_0))}$ is bounded
(cf. \cite{E1} and Remark 1.2 (iii)).
\end{remark}

To prove Theorem 3.1, we only need to show:

\begin{theorem}
Let $0<\alpha<1$ and $q>\frac{n}{1-\alpha}$. Suppose $Q_1\subset\Omega_{T}$, $(H)$ holds, and $u\in W^{2,1}_2(\Omega_T)$,
$Dp\in L_2(\Omega_T)$, $f\in L_2(\Omega_T)$ and $g\in L_\infty(-T,0;W^{1}_q(\Omega))$ satisfying (\ref{s0.0}) and
\begin{equation}\label{sh.24}
\sup_{t\in(-1,0]}\fint_{B_1}|u|^2dx\leq1~\mbox{and}~
\sup_{t\in(-1,0]}\fint_{B_1}|Dg|^qdx\leq1.
\end{equation}
There exist three constants $0<r_0<1$, $0<\varepsilon_0<1$ and $C>0$ depending only on $n$, $\lambda$, $\alpha$ and $q$, and three functions $a(t)\in L_\infty(-r^2_0,0;\mathbb{R}^{n})$, $b(t)\in L_\infty(-r^2_0,0;\mathbb{R}^{n^2})$ and $c(t)\in L_\infty(-r^2_0,0;\mathbb{R}^{n^3})$ with
\begin{equation}\label{j01}
\sup_{t\in(-r_0^{2},0]}|a(t)|,~~\sup_{t\in(-r_0^{2},0]}|b(t)|~~\mbox{and}\sup_{t\in(-r_0^{2},0]}|c(t)|\leq C
\end{equation}
such that if for any $0<r\leq1$,
\begin{equation}\label{st.21}
\fint_{Q_1}|D\nabla\times u|^2dxdt\leq r^2_0,
\end{equation}
\begin{equation}\label{sh.23}
\big|A(x,t)-A(0,t)\big|\leq\varepsilon_0r^{\alpha}~~\mbox{in}~~Q_r,
\end{equation}
\begin{equation}\label{n.25}
\begin{aligned}
\fint_{Q_r}|f-\overline{(f)}_{B_r}|^2dxdt\leq\varepsilon^2_0r^{2\alpha} ~\mbox{and}
\end{aligned}
\end{equation}
\begin{equation}\label{sh.25}
\begin{aligned}
\sup_{t\in(-r^2,0]}\fint_{B_r}|Dg-\overline{(Dg)}_{B_r}|^qdx\leq\varepsilon^q_0r^{q\alpha},
\end{aligned}
\end{equation}
then for any $0<r\leq r_0$,
\begin{equation}\label{sh.76}
\fint_{Q_{r}}|D\nabla\times u-\overline{(D\nabla\times u)}_{Q_{r}}|^2dxdt\leq Cr^{2\alpha}
\end{equation}
and
\begin{equation}\label{f4.09}
\begin{aligned}
\sup_{t\in(-r^{2},0]}\fint_{B_{r}}|u-a(t)-b(t)x-\frac12c(t)x^2|^2dx
\leq Cr^{2(2+\alpha)}.
\end{aligned}
\end{equation}
\end{theorem}

Theorem 3.1 follows from Theorem 3.3 by normalization. Actually, after suitable choice of coordinates, we may suppose $(x_0,t_0)=(0,0)$.
Choose $0<\tau<R$ such that $\tau^\alpha M_0\leq\varepsilon_0$.
For $(x,t)\in Q_1$, set
$$\tilde{A}(x,t)=A(\tau x,\tau^2t),~~\tilde{u}(x,t)=\frac{u(\tau x,\tau^2t)}{\tau^2 J},~~\tilde{p}(x,t)=\frac{p(\tau x,\tau^2t)}{\tau J},$$
$$\tilde{f}(x,t)=\frac{f(\tau x,\tau^2t)}{J}~~\mbox{and}~~\tilde{g}(x,t)=\frac{g(\tau x,\tau^2t)}{\tau J},$$
where
\begin{eqnarray*}
\begin{aligned}
J=&\ \frac{1}{\tau^2|B_\tau|^{\frac12}}\|u\|_{L_\infty(-T,0;L_2(\Omega))}
+\frac{1}{r_0}\frac{1}{|Q_\tau|^{\frac12}}\|D\nabla\times u\|_{L_2(\Omega_T)}\\
&+\frac{1}{|B_\tau|^{\frac1q}}\|g\|_{L_\infty(-T,0;W^1_q(\Omega))}
+\frac{1}{\varepsilon_0}M^{\frac12}_1+\frac{1}{\varepsilon_0}M^\frac{1}{q}_2.
\end{aligned}
\end{eqnarray*}
Then $\tilde u$, $\tilde p$, $\tilde f$ and $\tilde g$ satisfy (1.1) and for any $0<r\leq1$, 
\begin{equation*}
\big|\tilde{A}(x,t)-\tilde{A}(0,t)\big|=\big|A(\tau x,\tau^2t)-A(0,\tau^2t)\big|\leq M_0\tau^{\alpha}r^{\alpha}\leq\varepsilon_0r^{\alpha}~~\mbox{in}~~Q_r,
\end{equation*}
\begin{equation*}
\sup_{t\in(-1,0]}\fint_{B_1}|\tilde{u}|^2dx
=\frac{1}{\tau^{4}J^2}\sup_{t\in(-\tau^2,0]}\fint_{B_\tau}|u|^2dx\leq1,
\end{equation*}
\begin{equation*}
\fint_{Q_1}|D\nabla\times \tilde{u}|^2dxdt
=\frac{1}{J^2}\fint_{Q_\tau}|D\nabla\times u|^2dxdt\leq r_0^2,
\end{equation*}
\begin{equation*}
\begin{aligned}
\fint_{Q_r}|\tilde{f}-\overline{(\tilde{f})}_{B_r}|^2dxdt
=\frac{1}{J^2}\fint_{Q_{\tau r}}|f-\overline{(f)}_{B_{\tau r}}|^2dxdt\leq\varepsilon^2_0r^{2\alpha},
\end{aligned}
\end{equation*}
\begin{equation*}
\begin{aligned}
\sup_{t\in(-1,0]}\fint_{B_1}|D\tilde{g}|^qdx=\frac{1}{J^q}\sup_{t\in(-\tau^2,0]}\fint_{B_{\tau }}|Dg|^qdx\leq1
\end{aligned}
\end{equation*}
and
\begin{equation*}
\begin{aligned}
\sup_{t\in(-r^2,0]}\fint_{B_r}|D\tilde{g}-\overline{(D\tilde{g})}_{B_r}|^qdx
=\frac{1}{J^q}\sup_{t\in(-\tau^2r^2,0]}\fint_{B_{\tau r}}|Dg-\overline{(Dg)}_{B_{\tau r}}|^qdx\leq\varepsilon^q_0r^{q\alpha}.
\end{aligned}
\end{equation*}
By Theorem 3.3, there exist $\tilde{a}(t)\in L_\infty(-r^2_0,0;\mathbb{R}^{n})$, $\tilde{b}(t)\in L_\infty(-r^2_0,0;\mathbb{R}^{n^2})$ and $\tilde{c}(t)\in L_\infty(-r^2_0,0;\mathbb{R}^{n^3})$ such that
\begin{equation*}
\fint_{Q_{r}}|D\nabla\times \tilde{u}-\overline{(D\nabla\times \tilde{u})}_{Q_{r}}|^2dxdt\leq Cr^{2\alpha}
\end{equation*}
and
\begin{equation*}
\begin{aligned}
\sup_{t\in(-r^{2},0]}\fint_{B_{r}}|\tilde{u}-\tilde{a}(t)-\tilde{b}(t)x-\frac12\tilde{c}(t)x^2|^2dx
\leq Cr^{2(2+\alpha)}
\end{aligned}
\end{equation*}
for any $0<r<1$, where $C$ depends only on $n$, $\lambda$, $\alpha$ and $q$.
After rescaling back, we get (\ref{p01}) and (\ref{p02})
with
$a(t)=\tau^2 J\tilde{a}(\frac{t}{\tau^2})$, $b(t)=\tau J\tilde{b}(\frac{t}{\tau^2})$ and $c(t)=J\tilde{c}(\frac{t}{\tau^2})$.

\bigskip

Theorem 3.3 will be prove by the following iteration lemma.
\begin{lemma}
Under the hypotheses of Theorem 3.3, there exist three constants $0<r_0\leq\frac{1}{6}$, $0<\varepsilon_0<1$ and $C_0>0$ depending only on $n$, $\lambda$, $\alpha$ and $q$, and three functions $a_k(t)\in L_\infty(-r^{2k}_0,0;\mathbb{R}^{n})$, $b_k(t)\in L_\infty(-r^{2k}_0,0;\mathbb{R}^{n^2})$ and $c_k(t)\in L_\infty(-r^{2k}_0,0;\mathbb{R}^{n^3})$ for any $k\geq1$ such that
\begin{equation}\label{sh.27}
\sup_{t\in(-r_0^{2k},0]}|a_k(t)|\leq C_{a,k},~~~\sup_{t\in(-r_0^{2k},0]}|b_k(t)|\leq C_{b,k},~~~\sup_{t\in(-r_0^{2k},0]}|c_k(t)|
\leq C_{c,k},
\end{equation}
\begin{equation}\label{sh.26}
\fint_{Q_{r^k_0}}|D\nabla\times u-\overline{(D\nabla\times u)}_{Q_{r^k_0}}|^2dxdt\leq r_0^2r^{2k\alpha}_0,
\end{equation}
\begin{equation}\label{sh.28}
\begin{aligned}
\sup_{t\in(-r_0^{2k},0]}\fint_{B_{r^k_0}}|u-\overline{(u)}_{B_{r^k_0}}(t)-a_k(t)-b_k(t)x-\frac12c_k(t)x^2|^2dx
\leq r^{2k(2+\alpha)}_0,
\end{aligned}
\end{equation}
\begin{equation}\label{n0.11}
\begin{aligned}
\sup_{t\in(-r_0^{2k},0]}|\overline{(D\nabla\times u)}_{Q_{r^k_0}}-D\nabla\times\big(\frac12c_k(t)x^2\big)|
\leq r^{k\alpha}_0
\end{aligned}
\end{equation}
and
\begin{equation}\label{n0.12}
\begin{aligned}
\sup_{t\in(-r_0^{2k},0]}\fint_{B_{r^k_0}}|Dg-D\mbox{div}\big(\frac12c_k(t)x^2\big)|^qdx
\leq r^{qk\alpha}_0,
\end{aligned}
\end{equation}
where the constants $C_{a,k}$, $C_{b,k}$ and $C_{c,k}$ can be defined by the following way,
\begin{equation}\label{jj3}
C_{a,k}=r^{2k}_0C_{c,k},~~C_{b,k}=\sum_{j=1}^{k-1}\Big(C_{c,j}r_0^{1+j}+C_0r_0^{-\frac{n}2}r_0^{j(1+\alpha)}\Big)+C_0r_0^{-\frac{n}2}
\end{equation}
\begin{equation}\label{jj4}
\mbox{and}~~~C_{c,k}=C_0r^{-\frac n2}_0\sum^{k-1}_{j=0} r^{j\alpha}_0.
\end{equation}

\end{lemma}

\begin{remark}
(i) Observe that $C_{c,k}\leq C_0r_0^{-\frac n2}(1-r^\alpha_0)^{-1}$ and $C_{b,k}\leq2C_0r_0^{-\frac n2}(1-r^\alpha_0)^{-1}(1-r_0)^{-1}$ for any $k\geq1$. That is, $C_{a,k}$, $C_{b,k}$ and $C_{c,k}$ are all bounded.

(ii)  As mentioned before, since the coefficients of the second order polynomial $P_t(x):=\overline{(u)}_{B_{r^k_0}}(t)+a_k(t)+b_k(t)x+\frac12c_k(t)x^2$ may not be differentiable with respect to $t$, $u-P_t(x)$ will not satisfy the Stokes system and then is not a good object to be dealt with in the next scale as usual. We decompose the rescaled solution in the following Step 2 and 3 according to the structure of the equations, where to process the iteration, both an extra order of $r_0$ in (\ref{sh.26})
 and two additional estimates (\ref{n0.11}) and (\ref{n0.12}) are needed.
\end{remark}

\vskip 2mm\noindent{\bf Proof of Lemma 3.4.}\quad
The lemma is proved by induction on $k$. Since the proof of the conclusion as $k=1$ can be given by simplifying the proof of deriving the conclusion for $k+1$ from $k$, we postpone it
to the end of the proof.

The following is the proof of the conclusion for $k+1$ from $k$, which is divided into 7 steps and the constants $r_0$, $\varepsilon_0$ and $C_0$ will be determined in Step 6.

\emph{Step 0. Normalization.} For $(x,t)\in Q_1$, define
\begin{equation*}
\tilde{u}(x,t)=r_0^{-k(2+\alpha)}u(r^k_0x,r^{2k}_0t),
~~\tilde{p}(x,t)=r_0^{-k(1+\alpha)}p(r^k_0x,r^{2k}_0t)-r_0^{-k\alpha}\overline{(f)}_{B_{r^k_0}}(r^{2k}_0t)x,
\end{equation*}
\begin{equation*}
\tilde{A}(x,t)=A(r^k_0x,r^{2k}_0t),
\end{equation*}
\begin{equation*}
\tilde{f}(x,t)=r_0^{-k\alpha}\big(f(r^k_0x,r^{2k}_0t)-\overline{(f)}_{B_{r^k_0}}(r^{2k}_0t)\big)~~\mbox{and}~~\tilde{g}(x,t)=r_0^{-k(1+\alpha)}g(r^k_0x,r^{2k}_0t).
\end{equation*}
It follows that
\begin{eqnarray}\label{sh.58}
\left\{
\begin{array}{ll}
\vspace{2mm}
\tilde{u}_t-\tilde{A}(0,t)D^2\tilde{u}+\nabla \tilde{p}=\tilde{f}+\big(\tilde{A}(x,t)-\tilde{A}(0,t)\big)D^2\tilde{u}~~\mbox{in}~~Q_{1},\\
\vspace{2mm}
\mbox{div}~\tilde{u}=\tilde{g}~~\mbox{in}~~Q_{1}.
\end{array}
\right.
\end{eqnarray}

By (\ref{sh.23}), (\ref{n.25}) and (\ref{sh.25}),
\begin{equation}\label{sh.35}
\begin{aligned}
|\tilde{A}(x,t)-\tilde{A}(0,t)|=|A(r^k_0x,r^{2k}_0t)-A(0,r^{2k}_0t)|
\leq\varepsilon_0r_0^{k\alpha},
\end{aligned}
\end{equation}
\begin{eqnarray}\label{sh.40}
\begin{aligned}
\fint_{Q_1}|\tilde{f}|^2dxdt
=\fint_{Q_{r^k_0}}|r_0^{-k\alpha}\big(f(x,t)-\overline{(f)}_{B_{r^k_0}}(t)\big)|^2dxdt
\leq\varepsilon^2_0
\end{aligned}
\end{eqnarray}
and
\begin{eqnarray}\label{sh.90}
\begin{aligned}
\sup_{t\in(-r^2,0]}\fint_{B_r}|D\tilde{g}-\overline{(D\tilde{g})}_{B_r}|^qdx
&=r^{-kq\alpha}_0\sup_{t\in(-r^2r^{2k}_0,0]}\fint_{B_{rr^k_0}}|Dg(x,t)-\overline{(Dg)}_{B_{rr^k_0}}(t)|^qdx\\
&\leq\varepsilon^q_0r^{q\alpha}
\end{aligned}
\end{eqnarray}
for any $0<r\leq1$.

In view of the induction hypotheses (\ref{sh.26})-(\ref{n0.12}), we have
\begin{equation}\label{sh.30}
\fint_{Q_1}|D\nabla\times \tilde{u}-\overline{(D\nabla\times \tilde{u})}_{Q_1}|^2dxdt\leq r^2_0,
\end{equation}
\begin{equation}\label{sh.31}
\sup_{t\in(-1,0]}\fint_{B_1}|\tilde{u}-\alpha(t)-\beta(t)x-\frac12\gamma(t)x^2|^2dx\leq1,
\end{equation}
\begin{equation}\label{n0.13}
\begin{aligned}
\sup_{t\in(-1,0]}|\overline{(D\nabla\times \tilde{u})}_{Q_1}-D\nabla\times\big(\frac12\gamma(t)x^2\big)|
\leq 1
\end{aligned}
\end{equation}
and
\begin{equation}\label{n0.14}
\begin{aligned}
\sup_{t\in(-1,0]}\fint_{B_1}|D\tilde{g}-D\mbox{div}\big(\frac12\gamma(t)x^2\big)|^qdx
\leq 1,
\end{aligned}
\end{equation}
where
\begin{equation}\label{nj01}
\begin{aligned}
&\alpha(t)=r^{-k(2+\alpha)}_0\Big(\overline{(u)}_{B_{r^k_0}}(r^{2k}_0t)+a_k(r^{2k}_0t)\Big),
~~\beta(t)=r^{-k(1+\alpha)}_0b_k(r^{2k}_0t)\\
&\mbox{}\hskip3cm~~\mbox{and}~~\gamma(t)=r^{-k\alpha}_0c_k(r^{2k}_0t).
\end{aligned}
\end{equation}

\bigskip

\emph{Step 1. Estimates of $D^2\tilde{u}$.}~~
Applying Lemma 2.4 to $\tilde{u}-\alpha(t)-\beta(t)x-\frac12\gamma(t)x^2$
and by $\mbox{div}~\tilde{u}=\tilde{g}$,
\begin{eqnarray*}
\begin{aligned}
\lefteqn{\int_{Q_{\frac12}}|D^2\tilde{u}-\gamma(t)|^{2}dxdt= \int_{Q_{\frac12}}|D^2\big(\tilde{u}-\alpha(t)-\beta(t)x-\frac12\gamma(t)x^2\big)|^{2}}\hspace*{3mm}\\
\leq&\ C\Big(\int_{Q_1}|\nabla\times \tilde{u}-\overline{(D\nabla\times \tilde{u})}_{Q_{1}}|^{2}
+\int_{Q_1}|\overline{(D\nabla\times \tilde{u})}_{Q_{1}}-D\nabla\times\big(\frac12\gamma(t)x^2\big)|^{2}\\
&\mbox{}\hskip2cm +\int_{Q_1}|D\tilde{g}-D \mbox{div}\big(\frac12\gamma(t)x^2\big)|^{2}+\int_{Q_1}|\tilde{u}-\alpha(t)-\beta(t)x-\frac12\gamma(t)x^2|^{2}\Big),
\end{aligned}
\end{eqnarray*}
where $C$ depends only on $n$.
From (\ref{sh.30})-(\ref{n0.14}),
it follows that
\begin{eqnarray*}
\begin{aligned}
\int_{Q_{\frac12}}|D^2\tilde{u}-\gamma(t)|^{2}dxdt
\leq C,
\end{aligned}
\end{eqnarray*}
which combining with (\ref{sh.27}) leads to \begin{eqnarray}\label{j002}
\begin{aligned}
\int_{Q_{\frac12}}|D^2\tilde{u}|^{2}dxdt
\leq&\ C\sup_{t\in(-\frac{1}{4},0]}|r^{-k\alpha}_0c_k(r^{2k}_0t)|^{2}+C
\leq CC^2_{c,k}r^{-2k\alpha}_0,
\end{aligned}
\end{eqnarray}
where $C$ depends only on $n$.

Set $\tilde{h}(x,t)=\tilde f+\big(\tilde{A}(x,t)-\tilde{A}(0,t)\big)D^2\tilde{u}$ and then (\ref{sh.58}) reads as
\begin{eqnarray}\label{j003}
\left\{
\begin{array}{ll}
\vspace{2mm}
\tilde{u}_t-\tilde{A}(0,t)D^2\tilde{u}+\nabla \tilde{p}=\tilde h~~\mbox{in}~~Q_{1},\\
\vspace{2mm}
\mbox{div}~\tilde{u}=\tilde{g}~~\mbox{in}~~Q_{1}.
\end{array}
\right.
\end{eqnarray}
In view of (\ref{sh.35}), (\ref{sh.40}) and (\ref{j002}), we get
\begin{eqnarray}\label{sh.36}
\begin{aligned}
\int_{Q_{\frac12}}|\tilde{h}|^{2}dxdt
\leq&\ 2\int_{Q_{\frac12}}|\tilde f|^{2}dxdt+
2\int_{Q_{\frac12}}|\big(\tilde{A}(x,t)-\tilde{A}(0,t)\big)D^2\tilde{u}|^{2}dxdt\\
\leq &\ C\varepsilon_0^2+C\varepsilon_0^2r^{2k\alpha}_0C^2_{c,k}r^{-2k\alpha}_0
\leq C\varepsilon_0^2C^2_{c,k},
\end{aligned}
\end{eqnarray}
where $C$ depends only on $n$.

\bigskip

\emph{Step 2. Proof of (\ref{sh.26}) for $k+1$.} We decompose solution of (\ref{j003}) such that
$(\tilde{u},\tilde p)=(w, p_{w})+(v, p_{v})$
with
\begin{eqnarray}\label{sh.38}
\left\{
\begin{array}{ll}
\vspace{2mm}
w_t-\tilde{A}(0,t)D^2w+\nabla p_{w}=\tilde{h}~~\mbox{in}~~Q_{\frac12},\\
\vspace{2mm}
\mbox{div}~w=0~~\mbox{in}~~Q_{\frac12},\\
\vspace{2mm}
w=0~~\mbox{on}~~\partial_pQ_{\frac12}
\end{array}
\right.
\end{eqnarray}
and
\begin{eqnarray}\label{sh.39}
\left\{
\begin{array}{ll}
\vspace{2mm}
v_t-\tilde{A}(0,t)D^2v+\nabla p_{v}=0~~\mbox{in}~~Q_{\frac12},\\
\vspace{2mm}
\mbox{div}~v=\tilde{g}~~\mbox{in}~~Q_{\frac12},\\
\vspace{2mm}
v=\tilde{u}~~\mbox{on}~~\partial_pQ_{\frac12}.
\end{array}
\right.
\end{eqnarray}
Actually, we first solve $(w,p_w)$ from (\ref{sh.38}) by Lemma 2.6 and then set $(v, p_{v})=(\tilde{u},\tilde p)-(w, p_{w})$ which satisfies (\ref{sh.39}) clearly.

By $L_2$ estimates of (\ref{sh.38}),
\begin{eqnarray}\label{sh.45}
\begin{aligned}
\sup_{t\in(-\frac14,0]}\int_{B_{\frac12}}|w|^2dx+\int_{Q_{\frac12}}|Dw|^2dxdt
\leq C\int_{Q_{\frac12}}|\tilde{h}|^2dxdt
\leq C\varepsilon_0^2C^2_{c,k},
\end{aligned}
\end{eqnarray}
where (\ref{sh.36}) is used to derive the last inequality and $C$ depends only on $n$ and $\lambda$.

Taking $\nabla\times$ on both sides of the first equation in (\ref{sh.38}) and from
$\nabla\times\nabla p_{w}=0$, we deduce
\begin{eqnarray}\label{n0.22}
\begin{array}{ll}
\vspace{2mm}
\big(\nabla\times w\big)_t-\tilde{A}(0,t)D^2 \big(\nabla\times w\big)=\nabla\times\tilde{h}~~\mbox{in}~~Q_{\frac12}.
\end{array}
\end{eqnarray}
By Lemma 2.5 (ii), (\ref{sh.36}) and (\ref{sh.45}),
\begin{eqnarray}\label{n0.21}
\begin{aligned}
\|D\nabla\times w\|_{L_2(Q_{\frac13})}
\leq C\Big(\|\tilde{h}\|_{L_2(Q_{\frac12})}+\|\nabla\times w\|_{L_2(Q_{\frac12})}\Big)\leq C\varepsilon_0C_{c,k}.
\end{aligned}
\end{eqnarray}

Take $D\nabla\times$ on both sides of the first equation in (\ref{sh.39}) and denote
$$\tilde{v}=D\nabla\times v-\overline{(D\nabla\times v)}_{Q_{\frac13}}.$$ In view of (\ref{sh.30}), (\ref{n0.21}) and $v=\tilde{u}-w$,
\begin{eqnarray*}
\begin{aligned}
&\|\tilde{v}\|_{L_2(Q_{\frac13})}
 \leq \|D\nabla\times \tilde{u}-\overline{(D\nabla\times \tilde{u})}_{Q_{\frac13}}\|_{L_2(Q_{\frac13})}+\|D\nabla\times w-\overline{(D\nabla\times w)}_{Q_{\frac13}}\|_{L_2(Q_{\frac13})}\\
&\mbox{}\hskip2.5cm \leq C \Big(r_0+\varepsilon_0C_{c,k}\Big).
\end{aligned}
\end{eqnarray*}
By Lemma 2.5 (i), we have
\begin{eqnarray}\label{sh.j01}
\begin{aligned}
&\|\tilde{v}\|_{L_\infty(Q_{\frac14})}+\|D\tilde{v}\|_{L_\infty(Q_{\frac14})}
+\|\tilde{v}_t\|_{L_\infty(Q_{\frac14})}
\leq  C\|\tilde{v}\|_{L_2(Q_{\frac13})}\leq C \Big(r_0+\varepsilon_0C_{c,k}\Big),
\end{aligned}
\end{eqnarray}
where $C$ depends only on $n$ and $\lambda$.
Therefore
\begin{eqnarray*}
\begin{aligned}
\fint_{Q_{r_0}}|D\nabla\times v-\overline{(D\nabla\times v)}_{Q_{r_0}}|^2dxdt
& \leq Cr^2_0\Big(\|D\tilde{v}\|_{L_\infty(Q_{r_0})}
+\|\tilde{v}_t\|_{L_\infty(Q_{r_0})}\Big)^2\\
& \leq Cr^2_0\Big(r^2_0+\varepsilon_0^2C^2_{c,k}\Big).
\end{aligned}
\end{eqnarray*}

Combining it with (\ref{n0.21}),
\begin{equation*}
\begin{aligned}
&\fint_{Q_{r_0}}|D\nabla\times \tilde{u}-\overline{(D\nabla\times \tilde{u})}_{Q_{r_0}}|^2dxdt\\
&\leq 2\fint_{Q_{r_0}}|D\nabla\times w-\overline{(D\nabla\times w)}_{Q_{r_0}}|^2+2\fint_{Q_{r_0}}|D\nabla\times v-\overline{ (D\nabla\times v)}_{Q_{r_0}}|^2\\
&\leq C_1\varepsilon_0^2C^2_{c,k}r^{-(n+2)}_0
+C_2r^4_0
\end{aligned}
\end{equation*}
and
rescaling back,
\begin{eqnarray}\label{sh1.8}
\begin{aligned}
\fint_{Q_{r^{k+1}_0}}|D\nabla\times u-\overline{(D\nabla\times u)}_{Q_{r^{k+1}_0}}|^2dxdt\leq \Big(C_1\varepsilon_0^2C^2_{c,k}r^{-(n+2)}_0
+C_2r^4_0\Big)r_0^{2k\alpha},
\end{aligned}
\end{eqnarray}
where $C_1$ and $C_2$ depend only on $n$ and $\lambda$, which will give (\ref{sh.26}) for $k+1$ after all the constants are determined in Step 6.

\bigskip

\emph{Step 3. Proof of (\ref{sh.28}) for $k+1$.}
The inequalities (\ref{sh.31}), (\ref{sh.45}) and $v=\tilde{u}-w$ imply
\begin{eqnarray}\label{sh.47}
\begin{aligned}
 \lefteqn{\sup_{t\in(-\frac14,0]}\|v-\alpha(t)-\beta(t)x-\frac12\gamma(t)x^2\|_{L_2(B_{\frac12})}}\hspace*{14mm}\\
\leq&\ \sup_{t\in(-\frac14,0]}\|\tilde{u}-\alpha(t)-\beta(t)x-\frac12\gamma(t)x^2\|_{L_2(B_{\frac12})}
+\sup_{t\in(-\frac14,0]}\|w\|_{L_2(B_{\frac12})}\\
\leq&\ C\big(\varepsilon_0C_{c,k}+1\big),
\end{aligned}
\end{eqnarray}
where $C$ depends only on $n$ and $\lambda$.
By (\ref{sh.30}), (\ref{n0.13}) and (\ref{n0.21}),
\begin{eqnarray*}
\begin{aligned}
&\sup_{t\in(-\frac{1}{16},0]}|\overline{(D \nabla\times v)}_{Q_{\frac13}}-D\nabla\times (\frac12\gamma(t)x^2)|\\
&\mbox{}\hskip2cm \leq C\Big(\sup_{t\in(-\frac{1}{16},0]}|\overline{(D \nabla\times \tilde{u})}_{Q_{1}}-D\nabla\times (\frac12\gamma(t)x^2)|
+\|D \nabla\times w\|_{L_2(Q_{\frac13})}\\
&\mbox{}\hskip3cm +\|D \nabla\times \tilde{u}-\overline{(D \nabla\times \tilde{u})}_{Q_{1}}\|_{L_2(Q_{\frac13})}\Big)\\
&\mbox{}\hskip2cm \leq C\big(\varepsilon_0C_{c,k}+1\big).
\end{aligned}
\end{eqnarray*}
It follows from (\ref{sh.j01}) that
\begin{eqnarray*}
\begin{aligned}
&\sup_{t\in(-\frac{1}{16},0]}\|D\nabla\times v-D\nabla\times (\frac12\gamma(t)x^2)\|_{L_2(B_{\frac14})}\\
&\mbox{}\hskip2cm \leq
\sup_{t\in(-\frac{1}{16},0]}\Big(\|D\nabla\times v-\overline{(D \nabla\times v)}_{Q_{\frac13}}\|_{L_2(B_{\frac14})}\\
&\mbox{}\hskip3cm +\|\overline{(D \nabla\times v)}_{Q_{\frac13}}-D\nabla\times (\frac12\gamma(t)x^2)\|_{L_2(B_{\frac14})}\Big)\\
&\mbox{}\hskip2cm \leq C\big(\varepsilon_0C_{c,k}+1\big).
\end{aligned}
\end{eqnarray*}
Now, we apply Lemma 2.4 to $v-\alpha(t)-\beta(t)x-\frac12\gamma(t)x^2$ 
and then
\begin{eqnarray}\label{st.7}
\begin{aligned}
&\sup_{t\in(-\frac{1}{25},0]}\|D^2v-\gamma(t)\|_{L_2(B_{\frac15})}
 \leq C\sup_{t\in(-\frac{1}{16},0]}\Big(\|\mbox{div}\nabla\times v-\mbox{div}\nabla\times (\frac12\gamma(t)x^2)\|_{L_2(B_{\frac14})}\\
&\mbox{}\hskip0.8cm +\|D\tilde{g}-D\mbox{div}(\frac12\gamma(t)x^2)\|_{L_2(B_{\frac14})}
+\|v-\alpha(t)-\beta(t)x-\frac12\gamma(t)x^2\|_{L_2(B_{\frac14})}\Big)\\
&\mbox{}\hskip5cm \leq C\big(\varepsilon_0C_{c,k}+1\big),
\end{aligned}
\end{eqnarray}
where (\ref{n0.14}) and (\ref{sh.47}) are used to derive the last inequality and $C$ depends only on $n$ and $\lambda$.

Let
\begin{eqnarray*}
\begin{aligned}
\hat{v}=v-\alpha(t)-\beta(t)x-\frac12\overline{(D^2v)}_{B_{\frac15}}(t)x^2.
\end{aligned}
\end{eqnarray*}
Recall $\mbox{div} ~v=\tilde{g}$ and then
\begin{eqnarray*}
\begin{aligned}
\Delta \hat{v}=\mbox{div}\nabla\times v-\overline{(\mbox{div}\nabla\times v)}_{B_{\frac15}}+D \tilde{g}-\overline{(D \tilde{g})}_{B_{\frac15}}~~\mbox{in}~~B_{\frac15}.
\end{aligned}
\end{eqnarray*}
Set $\hat{v}=\hat{v}_1+\hat{v}_2$, where
\begin{eqnarray}\label{st.3}
\left\{
\begin{array}{ll}
\vspace{2mm}
\Delta \hat{v}_1=\mbox{div}\nabla\times v-\overline{(\mbox{div}\nabla\times v)}_{B_{\frac15}}+D \tilde{g}-\overline{(D \tilde{g})}_{B_{\frac15}}~~\mbox{in}~~B_{\frac15},\\
\vspace{2mm}
\hat{v}_1=0~~\mbox{on}~~\partial B_{\frac15}
\end{array}
\right.
\end{eqnarray}
and
\begin{eqnarray}\label{st.1}
\left\{
\begin{array}{ll}
\vspace{2mm}
\Delta \hat{v}_2=0~~\mbox{in}~~B_{\frac15},\\
\vspace{2mm}
\hat{v}_2=\hat{v}~~\mbox{on}~~\partial B_{\frac15}.
\end{array}
\right.
\end{eqnarray}

Since, by (\ref{sh.j01}),
\begin{eqnarray*}
\begin{aligned}
&\sup_{t\in(-\frac{1}{25},0]}\|\mbox{div}\nabla\times v-\overline{(\mbox{div}\nabla\times v)}_{B_{\frac15}}\|_{L_q(B_{\frac15})} \\
&\leq C\sup_{t\in(-\frac{1}{25},0]}\Big(\|D\nabla\times v-\overline{(D\nabla\times v)}_{Q_{\frac13}}\|_{L_q(B_{\frac15})} +\|\overline{(D\nabla\times v)}_{B_{\frac15}}-\overline{(D\nabla\times v)}_{Q_{\frac13}}\|_{L_q(B_{\frac15})}\Big)\\
&\leq C\|D\nabla\times v-\overline{(D\nabla\times v)}_{Q_{\frac13}}\|_{L_\infty(Q_{\frac15})}
\leq C\big(r_0+\varepsilon_0C_{c,k}\big),
\end{aligned}
\end{eqnarray*}
we have, by $W^2_q$ estimates for (\ref{st.3}) and (\ref{sh.90}),
\begin{eqnarray*}
\begin{aligned}
\sup_{t\in(-\frac{1}{25},0]}\|D^2\hat{v}_1\|_{L_q(B_{\frac15})}&\leq C\sup_{t\in(-\frac{1}{25},0]}\Big(\|\mbox{div}\nabla\times v-\overline{(\mbox{div}\nabla\times v)}_{B_{\frac15}}\|_{L_q(B_{\frac15})} \\
&\mbox{}\hskip2.8cm +\|D \tilde{g}-\overline{(D \tilde{g})}_{B_{\frac15}}\|_{L_q(B_{\frac15})}\Big)\\
&\leq C\big(r_0+\varepsilon_0C_{c,k}\big).
\end{aligned}
\end{eqnarray*}
It follows that
\begin{eqnarray}\label{st.10}
\begin{aligned}
\sup_{t\in(-r_0^2,0]}\fint_{B_{r_0}}|D^2\hat{v}_1-\overline{(D^2\hat{v}_1)}_{B_{r_0}}|^2dx
&\leq C\sup_{t\in(-r_0^2,0]} \fint_{B_{r_0}}|D^2\hat{v}_1|^2dx\\
&\leq C\sup_{t\in(-r_0^2,0]}r^{-\frac{2n}{q}}_0 \Big(\int_{B_{r_0}}|D^2\hat{v}_1|^qdx\Big)^{\frac{2}{q}}\\
&\leq Cr^{-\frac{2n}{q}}_0\big(r^2_0+\varepsilon^2_0C^2_{c,k}\big),
\end{aligned}
\end{eqnarray}
where $C$ depends only on $n$, $\lambda$ and $q$.

Since, by (\ref{sh.47}) and (\ref{st.7}),
\begin{eqnarray*}
\begin{aligned}
\sup_{t\in(-\frac{1}{25},0]}
\|\hat{v}\|_{L_2(B_{\frac15})}
&\leq \sup_{t\in(-\frac{1}{25},0]}\|v-\alpha(t)-\beta(t)x-\frac12\gamma(t)x^2\|_{L_2(B_{\frac15})}\\
&\mbox{}\hskip1cm +\sup_{t\in(-\frac{1}{25},0]}\|\frac12\overline{(D^2v)}_{B_{\frac15}}(t)x^2-\frac12\gamma(t)x^2\|_{L_2(B_{\frac15})}\\
&\leq C\big(\varepsilon_0C_{c,k}+1\big),
\end{aligned}
\end{eqnarray*}
we have, by the properties of harmonic functions and (\ref{st.1}),
\begin{eqnarray*}
\begin{aligned}
\sup_{t\in(-\frac{1}{36},0]}\|D^3\hat{v}_2\|_{L_\infty(B_{\frac16})}
\leq C\sup_{t\in(-\frac{1}{25},0]}
\|\hat{v}\|_{L_2(B_{\frac15})}
\leq C\big(\varepsilon_0C_{c,k}+1\big).
\end{aligned}
\end{eqnarray*}
It follows from Poincar\'{e}'s inequality that
\begin{eqnarray}\label{st.8}
\begin{aligned}
\sup_{t\in(-r_0^2,0]}\fint_{B_{r_0}}|D^2\hat{v}_2-\overline{(D^2\hat{v}_2)}_{B_{r_0}}|^2dx
&\leq C r_0^2\sup_{t\in(-r_0^2,0]}\|D^3\hat{v}_2\|^2_{L_\infty(B_{r_0})}\\
&\leq Cr_0^2\Big(\varepsilon^2_0C^2_{c,k}+1\Big).
\end{aligned}
\end{eqnarray}

Combining (\ref{st.8}) and (\ref{st.10}), we obtain
\begin{eqnarray}\label{st.5}
\begin{aligned}
&\sup_{t\in(-r_0^2,0]}\fint_{B_{r_0}}|D^2v-\overline{(D^2v)}_{B_{r_0}}(t)|^2dx
=\sup_{t\in(-r_0^2,0]}\fint_{B_{r_0}}|D^2\hat{v}-\overline{(D^2\hat{v})}_{B_{r_0}}(t)|^2dx\\
&\mbox{}\hskip0.1cm \leq \sup_{t\in(-r_0^2,0]} 2\Big(\fint_{B_{r_0}}|D^2\hat{v}_1-\overline{(D^2\hat{v}_1)}_{B_{r_0}}(t)|^2dx
+\fint_{B_{r_0}}|D^2\hat{v}_2-\overline{(D^2\hat{v}_2)}_{B_{r_0}}(t)|^2dx\Big)\\
&\mbox{}\hskip0.1cm
\leq Cr^{2-\frac{2n}{q}}_0+Cr^{-\frac{2n}{q}}_0\varepsilon^2_0C^2_{c,k},
\end{aligned}
\end{eqnarray}
where $C$ depends only on $n$, $\lambda$ and $q$.

Denote
\begin{eqnarray*}
\xi(t)=\fint_{B_{r_0}}\frac12\overline{(D^2v)}_{B_{r_0}}(t)x^2dx=\frac{nr^2_0}{2(n+2)}\overline{(D^2v)}_{B_{r_0}}(t).
\end{eqnarray*}
Since
\begin{eqnarray*}
\sup_{t\in(-r_0^2,0]}\fint_{B_{r_0}}\Big(v-\overline{(v)}_{B_{r_0}}(t)+\xi(t)-\overline{(D v)}_{B_{r_0}}(t)x-\frac12\overline{(D^2v)}_{B_{r_0}}(t)x^2\Big)dx=0
\end{eqnarray*}
and
\begin{eqnarray*}
\sup_{t\in(-r_0^2,0]}\fint_{B_{r_0}}\Big(Dv-\overline{(Dv)}_{B_{r_0}}(t)
-\overline{(D^2v)}_{B_{r_0}}(t)x\Big)dx=0,
\end{eqnarray*}
we have, by using Poincar\'{e}'s inequality twice and (\ref{st.5}),
\begin{eqnarray*}
\begin{aligned}
\lefteqn{\sup_{t\in(-r_0^2,0]}\fint_{B_{r_0}}|v-\overline{(v)}_{B_{r_0}}(t)+\xi(t)-\overline{(D v)}_{B_{r_0}}(t)x-\frac12\overline{(D^2v)}_{B_{r_0}}(t)x^2|^2dx}\hspace*{16mm}\\
\leq&\ Cr^2_0\sup_{t\in(-r_0^2,0]}\fint_{B_{r_0}}|Dv-\overline{(Dv)}_{B_{r_0}}(t)
-\overline{(D^2v)}_{B_{r_0}}(t)x|^2dx\\
\leq&\ Cr^4_0\sup_{t\in(-r_0^2,0]}\fint_{B_{r_0}}|D^2v-\overline{(D^2v)}_{B_{r_0}}(t)|^2dx\\
\leq&\ Cr^4_0\Big(r^{2-\frac{2n}{q}}_0+r^{-\frac{2n}{q}}_0\varepsilon^2_0C^2_{c,k}\Big).
\end{aligned}
\end{eqnarray*}

From
\begin{eqnarray*}
\begin{aligned}
\sup_{t\in(-r_0^2,0]}\fint_{B_{r_0}}|w-\overline{(w)}_{B_{r_0}}(t)|^2dx
\leq \frac{C}{r^{n}_0}\sup_{t\in(-\frac14,0]}\int_{B_{\frac12}}|w|^2dx
\leq C\varepsilon^2_0C^2_{c,k}r^{-n}_0
\end{aligned}
\end{eqnarray*}
implied by (\ref{sh.45}), it follows that
\begin{eqnarray*}
\begin{aligned}
&\sup_{t\in(-r_0^2,0]}\fint_{B_{r_0}}|\tilde{u}-\overline{(\tilde{u})}_{B_{r_0}}(t)+\xi(t)-\overline{(D v)}_{B_{r_0}}(t)x-\frac12\overline{(D^2v)}_{B_{r_0}}(t)x^2|^2dx\\
&\mbox{}\hskip1cm\leq 2\sup_{t\in(-r_0^2,0]}\left(\fint_{B_{r_0}}|v-\overline{(v)}_{B_{r_0}}(t)+\xi(t)-\overline{(D v)}_{B_{r_0}}(t)x-\frac12\overline{(D^2v)}_{B_{r_0}}(t)x^2|^2dx\right.\\
&\mbox{}\hskip4.5cm \left. +\fint_{B_{r_0}}|w-\overline{(w)}_{B_{r_0}}(t)|^2dx\right)\\
&\mbox{}\hskip1cm\leq C_3r^{6-\frac{2n}{q}}_0+C_4\varepsilon^2_0C^2_{c,k}r^{-n}_0,
\end{aligned}
\end{eqnarray*}
where $C_3$ and $C_4$ depend only on $n$, $\lambda$ and $q$.
Rescaling back,
\begin{equation}\label{ns8.2}
\begin{aligned}
&\sup_{t\in(-r_0^{2(k+1)},0]}\fint_{B_{r^{k+1}_0}}|u-\overline{(u)}_{B_{r^{k+1}_0}}(t)
-a_{k+1}(t)-b_{k+1}(t)x-\frac12c_{k+1}(t)x^2|^2dx\\
&\mbox{}\hskip2cm \leq \Big(C_3r^{6-\frac{2n}{q}}_0+C_4\varepsilon^2_0C^2_{c,k}r^{-n}_0\Big)r^{2k(2+\alpha)}_0
\end{aligned}
\end{equation}
with
\begin{eqnarray}\label{n28.2}
\begin{aligned}
&\mbox{}\hskip0.8cm a_{k+1}(t)=-r^{k(2+\alpha)}_0\xi\Big(\frac{t}{r^{2k}_0}\Big)
=-\frac{n}{2(n+2)}r^{k(2+\alpha)+2}_0\overline{(D^2v)}_{B_{r_0}}\Big(\frac{t}{r^{2k}_0}\Big),\\
&\mbox{}\hskip0.2cm~~b_{k+1}(t)=r^{k(1+\alpha)}_0\overline{(Dv)}_{B_{r_0}}\Big(\frac{t}{r^{2k}_0}\Big)
~~\mbox{and}~~c_{k+1}(t)=r^{k\alpha}_0\overline{(D^2v)}_{B_{r_0}}\Big(\frac{t}{r^{2k}_0}\Big).
\end{aligned}
\end{eqnarray}

Inequality (\ref{ns8.2}) will imply (\ref{sh.28}) for $k+1$ after all the constants are determined in Step 6.

\bigskip

\emph{Step 4. Proof of (\ref{sh.27}) for $k+1$.}~
In view of (\ref{sh.27}) and  (\ref{st.7}),
\begin{eqnarray*}
\begin{aligned}
&\sup_{t\in(-r_0^2,0]}\Big(\fint_{B_{r_0}}|D^2v|^2dx\Big)^{\frac12}
\leq \sup_{t\in(-r_0^2,0]}\left(\Big(\fint_{B_{r_0}}|\gamma(t)|^2dx\Big)^{\frac12}+\Big(\fint_{B_{r_0}}|D^2v-\gamma(t)|^2dx\Big)^{\frac12}\right)\\
&\mbox{}\hskip2cm\leq\sup_{t\in(-r_0^2,0]}|r^{-k\alpha}_0c_k(r^{2k}_0t)|
+\frac{C}{r^{\frac n2}_0}\sup_{t\in(-r_0^2,0]}\Big(\int_{B_{\frac15}}|D^2v-\gamma(t)|^2dx\Big)^{\frac12}\\
&\mbox{}\hskip2cm\leq C_{c,k}r^{-k\alpha}_0+C_5\big(\varepsilon_0C_{c,k}+1\big)r^{-\frac n2}_0,
\end{aligned}
\end{eqnarray*}
where $C_5$ depends only on $n$ and $\lambda$.
It follows that
\begin{eqnarray}\label{st.12}
\begin{aligned}
\sup_{t\in(-r_0^{2(k+1)},0]}|a_{k+1}(t)|
&\leq\frac{n}{2(n+2)}\sup_{t\in(-r_0^{2(k+1)},0]}r_0^{k(2+\alpha)+2}\Big(\fint_{B_{r_0}}|D^2 v(x,\frac{t}{r^{2k}_0})|^2dx\Big)^{\frac12}\\
&=\frac{n}{2(n+2)}\sup_{t\in(-r_0^2,0]}r_0^{k(2+\alpha)+2}\Big(\fint_{B_{r_0}}|D^2v|^2dx\Big)^{\frac12}\\
&\leq\frac{n}{2(n+2)}\Big(C_{c,k}+C_5\big(\varepsilon_0C_{c,k}+1\big)r^{-\frac n2}_0r_0^{k\alpha}\Big)r_0^{2k+2}
\end{aligned}
\end{eqnarray}
and
\begin{eqnarray}\label{n01.31}
\begin{aligned}
\sup_{t\in(-r_0^{2(k+1)},0]}|c_{k+1}(t)|
\leq C_{c,k}+C_5\big(\varepsilon_0C_{c,k}+1\big)r^{-\frac n2}_0r_0^{k\alpha}.
\end{aligned}
\end{eqnarray}

Since, by interpolation inequality,
\begin{eqnarray*}
\begin{aligned}
&\sup_{t\in(-\frac{1}{25},0]}\|Dv-\beta(t)-\gamma(t)x\|_{L_2(B_{\frac15})}\\
&\mbox{}\hskip1cm\leq C\sup_{t\in(-\frac{1}{25},0]}\Big(\|D^2v-\gamma(t)\|_{L_2(B_{\frac15})}
+\|v-\alpha(t)-\beta(t)x-\frac12\gamma(t)x^2\|_{L_2(B_{\frac15})}\Big),
\end{aligned}
\end{eqnarray*}
we have, by (\ref{sh.47}) and (\ref{st.7}),
\begin{eqnarray*}
\begin{aligned}
\sup_{t\in(-\frac{1}{25},0]}\|Dv-\beta(t)-\gamma(t)x\|_{L_2(B_{\frac15})}
\leq C\big(\varepsilon_0C_{c,k}+1\big).
\end{aligned}
\end{eqnarray*}
Combining with (\ref{sh.27}),
\begin{eqnarray*}
\begin{aligned}
\sup_{t\in(-r_0^2,0]}\Big(\fint_{B_{r_0}}|Dv|^2dx\Big)^{\frac12}
&\leq\sup_{t\in(-r_0^2,0]}\Big(\fint_{B_{r_0}}|\beta(t)|^2dx\Big)^{\frac12}
+\sup_{t\in(-r_0^2,0]}\Big(\fint_{B_{r_0}}|\gamma(t)x|^2dx\Big)^{\frac12}\\
&\ \mbox{}\hskip1cm+\sup_{t\in(-r_0^2,0]}\Big(\fint_{B_{r_0}}|Dv-\beta(t)-\gamma(t)x|^2dx\Big)^{\frac12}\\
&\leq C_{b,k}r^{-k(1+\alpha)}_0+C_{c,k}r^{1-k\alpha}_0+C_6\big(\varepsilon_0C_{c,k}+1\big)r^{-\frac n2}_0
\end{aligned}
\end{eqnarray*}
and then
\begin{eqnarray}\label{n02.31}
\begin{aligned}
\sup_{t\in(-r_0^{2(k+1)},0]}|b_{k+1}(t)|
&\leq \sup_{t\in(-r_0^{2(k+1)},0]}r_0^{k(1+\alpha)}\Big(\fint_{B_{r_0}}|Dv(x,\frac{t}{r^{2k}_0})|^2dx\Big)^{\frac12}\\
&\leq C_{b,k}+C_{c,k}r^{1+k}_0+C_6\big(\varepsilon_0C_{c,k}+1\big)r^{-\frac n2}_0r_0^{k(1+\alpha)},
\end{aligned}
\end{eqnarray}
where $C_6$ depends only on $n$ and $\lambda$.

Inequality (\ref{sh.27}) for $k+1$ will follow from (\ref{n01.31}) and (\ref{n02.31}) after all the constants are determined in Step 6.

\bigskip

\emph{Step 5. Proof of (\ref{n0.11}) and (\ref{n0.12}) for $k+1$.}
By (\ref{sh.j01}),
\begin{eqnarray*}
\begin{aligned}
&\|D\nabla\times v-\overline{(D\nabla\times v)}_{Q_{r_0}}\|_{L_\infty(Q_{r_0})}=\|\tilde{v}-\overline{(\tilde{v})}_{Q_{r_0}}\|_{L_\infty(Q_{r_0})}\\
&\mbox{}\hskip1cm\leq Cr_0\|D\tilde{v}\|_{L_\infty(Q_{\frac14})}+
Cr^2_0\|\tilde{v}_t\|_{L_\infty(Q_{\frac14})}\leq C\big(\varepsilon_0C_{c,k}+1\big)r_0
\end{aligned}
\end{eqnarray*}
and then by (\ref{n0.21}),
\begin{equation*}
\begin{aligned}
& \sup_{t\in(-r_0^2,0]}|\overline{(D\nabla\times v)}_{B_{r_0}}(t)-\overline{(D\nabla\times \tilde{u})}_{Q_{r_0}}|
\leq \|D\nabla\times v-\overline{(D\nabla\times \tilde{u})}_{Q_{r_0}}\|_{L_\infty(Q_{r_0})}
\\
&\mbox{}\hskip0.5cm \leq \|D\nabla\times v-\overline{(D\nabla\times v)}_{Q_{r_0}}\|_{L_\infty(Q_{r_0})}+
|\overline{(D\nabla\times w)}_{Q_{r_0}}|
\leq C_7r_0+ C_8\varepsilon_0C_{c,k}r^{-\frac{n+2}{2}}_0,
\end{aligned}
\end{equation*}
where $C_7$ and $C_8$ depend only on $n$ and $\lambda$.
Rescaling back and by (\ref{n28.2}),
\begin{equation}\label{ns2.9}
\begin{aligned}
&\sup_{t\in(-r_0^{2(k+1)},0]}|\overline{(D\nabla\times u)}_{Q_{r^{k+1}_0}}-D\nabla\times\big(\frac12c_{k+1}(t)x^2\big)|\\
&\mbox{}\hskip2cm \leq \Big(C_7r_0+ C_8\varepsilon_0C_{c,k}r^{-\frac{n+2}{2}}_0\Big)r^{k\alpha}_0.
\end{aligned}
\end{equation}

By (\ref{sh.90}),
\begin{equation*}
\begin{aligned}
\sup_{t\in(-r_0^2,0]}\fint_{B_{r_0}}|D\tilde{g}-\overline{(D\tilde{g})}_{B_{r_0}}(t)|^qdx
&\leq \varepsilon^q_0r^{q\alpha}_0.
\end{aligned}
\end{equation*}
Rescaling back and by (\ref{n28.2}) and $\mbox{div}~v=\tilde{g}$,
\begin{equation}\label{ns2.4}
\begin{aligned}
\sup_{t\in(-r_0^{2(k+1)},0]}\fint_{B_{r^{k+1}_0}}|Dg-D\mbox{div}\big(\frac12c_{k+1}(t)x^2\big)|^qdx
\leq \varepsilon^q_0r^{q(k+1)\alpha}_0.
\end{aligned}
\end{equation}

After all the constants are determined in Step 6, we will see (\ref{n0.11}) and (\ref{n0.12}) for $k+1$ via (\ref{ns2.9}) and (\ref{ns2.4}) respectively.
\bigskip

\emph{Step 6. Determination of the constants.}
Let $C_i$ ($i=1\sim 8$)  be given by
(\ref{sh1.8}), (\ref{ns8.2}), (\ref{n01.31}), (\ref{n02.31}) and (\ref{ns2.9})
respectively. We first set $$C_0=\max\{ 2C_5, 2C_6\}$$ and choose $0<r_0\leq\frac16$ small enough such that
\begin{eqnarray*}
C_2r^4_0\leq \frac12r_0^{2\alpha},~~C_3r^{6-\frac{2n}{q}}_0\leq\frac12r_0^{2(2+\alpha)}
~~\mbox{and}~~C_7r_0\leq \frac12r_0^{\alpha}.
\end{eqnarray*}
Next, we choose $0<\varepsilon_0<1$ small enough such that
\begin{equation}\label{n7.144}
\max\{C_1,C_4,C^2_8\}C_*^2r^{-(n+2)}_0\varepsilon_0^2\leq \frac12r_0^{2(2+\alpha)},	
\end{equation}
where $C_*=C_0r_0^{-\frac n2}(1-r^\alpha_0)^{-1}$.

By (\ref{n7.144}), we see (\ref{sh.26}), (\ref{sh.28}), (\ref{n0.11}) and (\ref{n0.12}) for $k+1$ clearly via (\ref{sh1.8}), (\ref{ns8.2}), (\ref{ns2.9}) and (\ref{ns2.4}) respectively.

Set $C_{b,1}=C_0r^{-\frac n2}_0$ and $C_{c,1}=C_0r^{-\frac n2}_0$.
Since $\varepsilon_0C_*\leq1$ (implied by (\ref{n7.144}) clearly), we deduce from (\ref{n01.31}) and (\ref{n02.31}),
\begin{eqnarray*}
\begin{aligned}
\sup_{t\in(-r_0^{2(k+1)},0]}|c_{k+1}(t)|
\leq C_{c,k}+C_5\big(\varepsilon_0C_*+1\big)r^{-\frac n2}_0r_0^{k\alpha}\leq C_{c,k}+C_0r^{-\frac n2}_0r_0^{k\alpha}=C_{c,k+1}
\end{aligned}
\end{eqnarray*}
and
\begin{eqnarray*}
\begin{aligned}
\sup_{t\in(-r_0^{2(k+1)},0]}|b_{k+1}(t)|
\leq&\ C_{b,k}+C_{c,k}r^{1+k}_0+C_6\big(\varepsilon_0C_*+1\big)r^{-\frac n2}_0r_0^{k(1+\alpha)}\\
\leq&\ C_{b,k}+C_{c,k}r^{1+k}_0+C_0r^{-\frac n2}_0r_0^{k(1+\alpha)}
=C_{b,k+1}.
\end{aligned}
\end{eqnarray*}
Since $a_{k+1}(t)=-\frac{n}{2(n+2)}r^{2k+2}_0c_{k+1}(t)$,
we see (\ref{sh.27}), (\ref{jj3}) and (\ref{jj4}) clearly.

\bigskip

\emph{Finally, we simplify the above arguments to show the conclusion as $k=1$.} Actually, we only need to repeat the above Step 1 to 5 by setting $\alpha,\beta,\gamma=0$ and
using (\ref{sh.24}) and (\ref{st.21}) instead of the induction assumptions on $k$.
\qed

\bigskip

\noindent{\bf Proof of Theorem 3.3.}\quad
For any $0<r\leq r_0$, let $k\geq1$ such that
\begin{equation}\label{j05}
r^{k+1}_0< r\leq r^k_0.
\end{equation}
By (\ref{sh.26}), we have
\begin{equation*}
\begin{aligned}
&\fint_{Q_{r}}|D\nabla\times u-\overline{(D\nabla\times u)}_{Q_{r}}|^2dxdt\\
&\mbox{}\hskip1cm\leq\frac{|Q_{r^k_0}|}{|Q_{r}|}\fint_{Q_{r^k_0}}|D\nabla\times u-\overline{(D\nabla\times u)}_{Q_{r^k_0}}|^2dxdt
\leq Cr^{2k\alpha}_0
\leq Cr^{2\alpha}
\end{aligned}
\end{equation*}
and hence (\ref{sh.76}) holds.

For any $t\in(-r^2_0,0)$, define
\begin{equation*}
a(t)=\overline{(u)}_{B_{r^k_0}}(t)+a_k(t),
\end{equation*}
\begin{equation*}
b(t)=b_k(t)~~ \mbox{and}~~c(t)=c_k(t)
\end{equation*}
for $t\in(-r^{2k}_0,-r^{2(k+1)}_0]$ and $k\geq 1$, where $a_k$, $b_k$ and $c_k$ are given by Lemma 3.4.

From (\ref{sh.28}) and Lemma 2.3,
\begin{equation*}
\begin{aligned}
\sup_{t\in(-r^{2(k+1)}_0,0]}|\overline{(u)}_{B_{r^{k}_0}}(t)+a_k(t)
-\overline{(u)}_{B_{r^{k+1}_0}}(t)-a_{k+1}(t)|
\leq&\ Cr^{k(2+\alpha)}_0,
\end{aligned}
\end{equation*}
\begin{equation*}
\sup_{t\in(-r^{2(k+1)}_0,0]}|b_{k}(t)-b_{k+1}(t)|\leq Cr^{k(1+\alpha)}_0
\end{equation*}
and
\begin{equation*}
\sup_{t\in(-r^{2(k+1)}_0,0]}|c_{k}(t)-c_{k+1}(t)|\leq Cr^{k\alpha}_0.
\end{equation*}
Define
$$a(0)=\lim_{k\rightarrow\infty}\overline{(u)}_{B_{r^{k}_0}}(0)+a_k(0),~~
b(0)=\lim_{k\rightarrow\infty}b_{k}(0)~~\mbox{and}~~c(0)=\lim_{k\rightarrow\infty}c_{k}(0).$$

Since, by the definition of $a(t)$,
\begin{equation*}
\begin{aligned}
\lefteqn{\sup_{t\in(-r^{2k}_0,0]}|a(t)-\overline{(u)}_{B_{r^k_0}}(t)-a_k(t)|
=\sup_{t\in(-r^{2(k+1)}_0,0]}|a(t)-\overline{(u)}_{B_{r^k_0}}(t)-a_k(t)|}\hspace*{10mm}\\
\leq&\ \sup_{t\in(-r^{2(k+1)}_0,0]}\Big(|a(t)-\overline{(u)}_{B_{r^{k+1}_0}}(t)-a_{k+1}(t)|\\
&\mbox{}\hskip2.5cm +
|\overline{(u)}_{B_{r^k_0}}(t)+a_k(t)-\overline{(u)}_{B_{r^{k+1}_0}}(t)-a_{k+1}(t)|
\Big)\\
\leq&\
\sup_{t\in(-r^{2(k+1)}_0,0]}|a(t)-\overline{(u)}_{B_{r^{k+1}_0}}(t)-a_{k+1}(t)|+
Cr^{k(2+\alpha)}_0,
\end{aligned}
\end{equation*}
we have, by repeating the above inequalities,
\begin{equation}\label{sh2.71}
\sup_{t\in(-r^{2k}_0,0]}|a(t)-\overline{(u)}_{B_{r^k_0}}(t)-a_k(t)|\leq C\sum_{j=k}^\infty r_0^{j(2+\alpha)}\leq Cr_0^{k(2+\alpha)}
\end{equation}
and similarly,
\begin{equation*}
\sup_{t\in(-r^{2k}_0,0]}|b(t)-b_{k}(t)|\leq C\sum_{j=k}^\infty r_0^{j(1+\alpha)}\leq Cr_0^{k(1+\alpha)}
\end{equation*}
and
\begin{equation}\label{n0.52}
\sup_{t\in(-r^{2k}_0,0]}|c(t)-c_{k}(t)|\leq C\sum_{j=k}^\infty r_0^{j\alpha}\leq Cr_0^{k\alpha},
\end{equation}
where $C$ depends only on $n$, $\lambda$, $\alpha$ and $q$.
By (\ref{sh.28}) again,
\begin{equation*}
\begin{aligned}
\sup_{t\in(-r^{2k}_0,0]}\fint_{B_{r_0^k}}|u-a(t)-b(t)x-\frac12c(t)x^2|^2dx
\leq Cr_0^{2k(2+\alpha)}.
\end{aligned}
\end{equation*}
For any $0<r\leq r_0$, let $k$ be given by (\ref{j05}) and it follows that
\begin{equation*}
\begin{aligned}
&\sup_{t\in(-r^{2},0]}\fint_{B_{r}}|u-a(t)-b(t)x-\frac12c(t)x^2|^2dx\\
&\mbox{}\hskip0.2cm \leq r_0^{-n}\sup_{t\in(-r^{2k}_0,0]}\fint_{B_{r_0^k}}|u-a(t)-b(t)x-\frac12c(t)x^2|^2dx
\leq Cr_0^{2k(2+\alpha)}\leq Cr^{2(2+\alpha)}.
\end{aligned}
\end{equation*}
Thus (\ref{f4.09}) is proved.

Now, we show the boundedness of $a(t)$, $b(t)$ and $c(t)$. Actually,
by (\ref{sh.27}) and Remark 3.5 (i), we obtain
\begin{equation*}
\sup_{t\in(-r^{2}_0,0]}|b(t)|\leq\sup_k\sup_{t\in(-r^{2k}_0,-r^{2(k+1)}_0]}|b_{k}(t)|
\leq\sup_k C_{b,k}
\leq C,
\end{equation*}
\begin{equation*}
\sup_{t\in(-r^{2}_0,0]}|c(t)|\leq\sup_k\sup_{t\in(-r^{2k}_0,-r^{2(k+1)}_0]}|c_{k}(t)|
\leq\sup_k C_{c,k}
\leq C,
\end{equation*}
and from (\ref{sh.24}) and (\ref{sh2.71}), we see
\begin{equation*}
\begin{aligned}
\sup_{t\in(-r^{2}_0,0]}|a(t)|=&\
\sup_{t\in(-r^{2}_0,0]}\Big(|\overline{(u)}_{B_{r_0}}(t)|+|a_1(t)|+
|a(t)-\overline{(u)}_{B_{r_0}}(t)-a_1(t)|\Big)
\\
\leq&\ Cr^{-\frac n2}_0\sup_{t\in(-r^{2}_0,0]}\Big(\int_{B_{1}}|u|^2dx\Big)^{\frac12}+Cr^{-\frac n2}_0+Cr_0^{2+\alpha}.
\end{aligned}
\end{equation*}
The proof of Theorem 3.3 is complete.
\qed

\bigskip

\noindent{\bf Proof of Theorem 1.1 (i).}\quad Let $R=\min\{\mbox{dist}(\Omega',\partial\Omega),\sqrt{T-T'}\}$ and $q=\frac{2n}{1-\alpha}$. Set
$M_1=[f]_{x,\alpha;\Omega_T}$ and $M_2=[Dg]_{x,\alpha;\Omega_T}$.

For any $(x_0,t_0)\in\Omega'_{T'}$, by Theorem 3.1, there exists
$\sigma$ depending only on $n$, $\lambda$, $\alpha$ and $M_0$ such that (\ref{p01}) holds for any $0<r\leq r_0=\sigma R$. From Lemma 2.1 (i), we obtain the estimate of $[D\nabla\times u]_{\alpha,\frac{\alpha}{2};\overline{\Omega'_{T'}}}$. Then the desired estimate
for $\|D\nabla\times u\|_{C^{\alpha,\frac{\alpha}{2}}(\overline{\Omega'_{T'}})}$ in (\ref{g01}) follows from interpolation and the boundedness of $\|D\nabla\times u\|_{L_2(\Omega_{T})}$ immediately.

By Theorem 3.1 again, there exist
$a(x_0,t)\in L_\infty(t_0-r_0^2,t_0;\mathbb{R}^{n})$, $b(x_0,t)\in L_\infty(t_0-r_0^2,t_0;\mathbb{R}^{n^2})$ and $c(x_0,t)\in L_\infty(t_0-r_0^2,t_0;\mathbb{R}^{n^3})$
such that (\ref{p02}) holds for any $0<r\leq r_0$. Set $t=t_0$ in the obtained (\ref{p02}) and then
\begin{equation*}
\begin{aligned}
\left(\fint_{B_{r}(x_0)}|u(x,t_0)-a(x_0,t_0)-b(x_0,t_0)x-\frac12c(x_0,t_0)x^2|^2dx\right)^{\frac12}
\leq&\ Cr^{2+\alpha}{\cal K}
\end{aligned}
\end{equation*}
for any $0<r\leq r_0$.
Then by Lemma 2.1 (ii) and interpolation, we obtain the desired estimate of $\|D^2u\|_{C_x^{\alpha}(\overline{\Omega'_{T'}})}$ in  (\ref{g01}).
\qed


\bigskip

\section{Estimates of pressure}\setcounter{equation}{0}
This section is devoted to first order H\"{o}lder estimates for the pressure $p$, that is, the proof of Theorem 1.1 (ii).

\bigskip

\noindent{\bf Proof of Theorem 1.1 (ii).}\quad
Suppose $Q_1\subset\Omega_{T}$. Taking $\mbox{div}$ on both sides of the first equation in (\ref{s0.0}) and by $\mbox{div}~ u=g$, we obtain
\begin{equation}\label{n5.14}
\Delta p=\mbox{div}\big(f+A(x,t)D^2u\big)-g_t~~\mbox{in}~~Q_1.
\end{equation}
Note that $A\in C_{x}^{\alpha}(Q_1)$, $f\in C_{x}^{\alpha}(Q_1)$, $p\in L_\infty(-1,0;L_2(B_1))$ and $g_t\in L_\infty(-1,0;L_{q}(B_1))$ and $D^2u\in C_{x}^{\alpha}(\overline{Q_{\frac12}})$ by Theorem 1.1 (i). Then, using $C^{1,\alpha}$ estimates for elliptic equations (\ref{n5.14}) and taking supremum with respect to $t$, we obtain $Dp\in C_{x}^{\alpha}(\overline{Q_{\frac14}})$ and
\begin{equation*}
\begin{aligned}
\|Dp\|_{C_x^{\alpha}(\overline{Q_{\frac14}})}
\leq
C\Big(\|D^2u\|_{C_{x}^{\alpha}(Q_{\frac12})}
+\|f\|_{C_{x}^{\alpha}(Q_{\frac12})}+\|g_t\|_{L_\infty(-\frac{1}{4},0;L_{q}(B_{\frac12}))}
+\|p\|_{L_\infty(-\frac{1}{4},0;L_2(B_{\frac12}))}\Big),
\end{aligned}
\end{equation*}
where $C$ depends only on $n$, $\lambda$, $\alpha$ and $M_0$.
By (\ref{g01}) and covering arguments, Theorem 1.1 (ii) follows immediately.
\qed

\begin{remark}

(i) If $g_t=\mbox{div}~G$ for some $G\in L_{\infty}(-T,0; C_x^{\alpha}(\Omega))$, then  from (\ref{n5.14}), we still have $Dp\in C_{x}^{\alpha}(\overline{Q_{\frac14}})$ and
\begin{equation*}
\begin{aligned}
\|Dp\|_{C_x^{\alpha}(\overline{Q_{\frac14}})}
\leq
C\Big(\|D^2u\|_{C_{x}^{\alpha}(Q_{\frac12})}
+\|f\|_{C_{x}^{\alpha}(Q_{\frac12})}+\|G\|_{C_{x}^{\alpha}(Q_{\frac12})}
+\|p\|_{L_\infty(-\frac{1}{4},0;L_2(B_{\frac12}))}\Big)
\end{aligned}
\end{equation*}
for some $C$ depending only on $n$, $\lambda$, $\alpha$ and $M_0$.

(ii) From the pointwise estimates (\ref{p01}) and (\ref{p02}), one can only derive
\begin{equation*}
\begin{aligned}
\fint_{B_{r}}\left|\int_{-T}^0\left(Dp-\overline{(Dp)}_{B_r}\right)dt\right|^2dx
\leq Cr^{2\alpha}
\end{aligned}
\end{equation*}
for some constant $C$ independent of $r$,
which is weaker than the pointwise $C^{\alpha}_x$ estimate. Observe (\ref{n5.14}) can not give regularity for $p$ in the time direction.
\end{remark}

\bigskip

\centerline{REFERENCES}
\bibliographystyle{plain}

\end{document}